\documentclass[final,3p,times]{modified-elsarticle}
\pdfoutput=1

\usepackage{amssymb,amsmath,amsfonts,mathtools}
\usepackage{amsthm}
	\theoremstyle{remark}
	\newtheorem{rem}{Remark}[section]

	\theoremstyle{plain}

	\DeclareMathAlphabet{\mathbcal}{OMS}{cmsy}{b}{n}
	
\usepackage{algorithm}
\usepackage{algorithmicx,algpseudocode}
\usepackage{cases}

\usepackage[nodots]{numcompress}

\usepackage{hyperref}
\usepackage[]{setspace}

\biboptions{sort&compress}
\usepackage[normalem]{ulem}

\journal{arXiv preprint}

\renewcommand{\vec}[1]{\mathbf{#1}}
\newcommand{\tensorTwo}[1]{\boldsymbol#1}
\newcommand{\tensorFour}[1]{\mathbb{#1}}

\newcommand{\Mat}[1]{\mathsf{#1}}
\renewcommand{\Vec}[1]{\mathsf{#1}}

\newcommand{\blkMat}[1]{\boldsymbol{\mathsf{#1}}}
\newcommand{\blkVec}[1]{\boldsymbol{\mathsf{#1}}}

\newcommand{\sub}[1]{_\mathit{#1}}

\usepackage{booktabs}
\usepackage{color, soul}
\usepackage{subcaption}
\usepackage{arydshln}
\usepackage{upgreek}
\usepackage{hhline}
\usepackage{multirow}

\begin{document}

\begin{frontmatter}



\title{A Two-Stage Preconditioner for Multiphase Poromechanics in Reservoir Simulation}


\author[label1]{Joshua A. White\corref{cor1}}
\ead{jawhite@llnl.gov}

\author[label1,label2]{Nicola Castelletto}
\ead{castelletto1@llnl.gov}

\author[label2]{Sergey Klevtsov}
\ead{klevtsov@stanford.edu}

\author[label3]{Quan M. Bui}
\ead{ bui9@llnl.gov}

\author[label3]{Daniel Osei-Kuffuor}
\ead{ oseikuffuor1@llnl.gov}

\author[label2]{Hamdi A. Tchelepi}
\ead{tchelepi@stanford.edu}

\address[label1]{Atmospheric, Earth and Energy Division,
                 Lawrence Livermore National Laboratory,
                 United States}

\address[label2]{Energy Resources Engineering,
                 Stanford University,
                 United States}

\address[label3]{Center for Applied Scientific Computing,
                 Lawrence Livermore National Laboratory,
                 United States}
                 
\cortext[cor1]{Corresponding author}

\begin{abstract}
Many applications involving porous media---notably reservoir engineering and geologic applications---involve tight coupling between multiphase fluid flow, transport, and poromechanical deformation.
While numerical models for these processes have become commonplace in research and industry, the poor scalability of existing solution algorithms has limited the size and resolution of models that may be practically solved.  In this work, we propose a two-stage Newton-Krylov solution algorithm to address this shortfall.  The proposed solver exhibits rapid convergence, good parallel scalability, and is robust in the presence of highly heterogeneous material properties.  The key to success of the solver is a block-preconditioning strategy that breaks the fully-coupled system of mass and momentum balance equations into simpler sub-problems that may be readily addressed using targeted algebraic methods.  Numerical results are presented to illustrate the performance of the solver on challenging benchmark problems.
\end{abstract}


\end{frontmatter}


\allowdisplaybreaks


\section{Introduction}
\label{sec:intro}

Coupled simulations of multiphase fluid flow, transport, and geomechanics are useful for predicting performance of many subsurface systems.   These simulations solve a set of nonlinear, time-dependent partial differential equations (PDEs) enforcing mass conservation of each fluid phase and linear momentum balance for the mixture.    For this system, fully-implicit time discretization is widely preferred for its unconditional stability, allowing large timesteps.  The major computational bottleneck, however, is that one must solve a sequence of large, ill-conditioned linear systems to advance the solution.   Robust and scalable solvers are therefore needed to enable large-scale simulations on high performance computing platforms.   In this paper, we propose an efficient preconditioning strategy for multiphase poromechanics that addresses the inherent ill-conditioning of the discrete system using a mixture of physics-based insight and flexible algebraic methods. 

Over the last four decades, a rich literature focusing on the development of efficient iterative solvers for fully-implicit simulation of complex multiphase flow and transport---without mechanics---has developed within the reservoir simulation community.   The most popular approach is based on Krylov subspace methods combined with a  \textit{Constrained Pressure Residual} (CPR) multistage preconditioning technique \cite{Wal83,WalKenLit85}.   In essence, the CPR preconditioner uses separate stages to damp error modes associated with (i) the nearly-hyperbolic behavior of the saturation (or composition) variables and (ii) the nearly-elliptic character of the pressure variables.   The standard reservoir simulator implementation \cite{LacVasWhe01,SchMasWen03,Cao_etal05,Stu_etal07,ZhoJiaThc13,Gri_etal14,LiuWanChe16} couples Incomplete LU (ILU) factorizations, effective for the hyperbolic part, with an algebraic multigrid (AMG) method, well-suited for the elliptic part.   A variant of CPR preconditioning that replaces AMG by a multiscale finite volume (MSFV) or finite element (MSFE) solver was proposed in \cite{Cus_etal15}.   Using a multigrid reduction (MGR) approach \cite{RieTroWin83}, the CPR method has been generalized in a broader multigrid framework and successfully applied to multiphase flow and transport with phase transitions in \cite{Wan_etal17,BuiWanOse18}.   Recent solution alternatives to CPR can be found in \cite{Yan_etal18,SinPenWhe18}.

To date, the design of fully-implicit solvers for poromechanics has primarily targeted single-phase flow in fully or partially saturated media.   Most efforts can be grouped in two broad  categories: (i) preconditioned Krylov solvers  \cite{BerFerGam07,WhiBor11,BerMar12,HagOsnLan12a,WhiCasTch16,LeeMarWin17,Adl_etal18}, and (ii) multigrid solvers \cite{Gas_etal04,Luo_etal15,GasRod17,Luo_etal17}---with obvious overlap.   In particular, the analysis of block preconditioners for the displacement-velocity-pressure formulation has been given increasing attention during  recent years \cite{Lip02,Kuz_etal03,FerCasGam10,TurArb14,CasWhiFer16}.   The extension to the case of multiphase flow and transport has typically been addressed by adopting sequential-implicit approaches \cite{SetMou98,ManLon02,Dai_etal02,Jea_etal_07,KimTchJua11a,MikWhe13,KimTchJua13,DosNor15,CorMur18} where easier sub-problems are solved separately but iterated until coupled convergence.   In this approach, a proper splitting strategy is necessary to guarantee unconditional convergence of the sequential iteration.  The most established algorithms rely on the \textit{fixed-stress} coupling scheme \cite{KimTchJua11a,KimTchJua11b,MikWhe13}.


In the present work, a preconditioning scheme for fully-implicit multiphase fluid flow, transport, and geomechanics is proposed.   The algorithm relies on successive sparse Schur-complement approximations to reduce the coupled system into easier sub-problems for which algebraic methods can be applied.  We build on fixed-stress partitioning \cite{SetMou98,KimTchJua11a,MikWhe13,WhiCasTch16} and CPR \citep{Wal83} concepts to construct a novel multi-stage preconditioning scheme well-suited to the problem at hand. The end result is a modular, easy-to-implement, and highly scalable solver. 
 
 The paper is organized as follows.    Sections~\ref{sec:model} and \ref{sec:num_model} introduce the mathematical model and discretization scheme.  The nonlinear and linear solution algorithms are described in detail in Sections~\ref{sec:NewKry_solver} and \ref{sec:precond}. Weak and strong scaling results are presented in Section~\ref{sec:numerical_results} to demonstrate performance and robustness of the proposed preconditioner in a parallel setting.  We then end with a few concluding remarks regarding future work.


\section{Governing Equations}
\label{sec:model}
As a simple but representative model system, we consider a displacement-saturation-pressure formulation of multiphase poroelasticity~\citep{Cou04} assuming isothermal, quasi-static, two-phase flow of immiscible fluids. We limit the discussion to small-strain kinematics.  Let $\overline{\Omega} = \Omega \cup \Gamma$ denote a closed domain in $\mathbb{R}^3$ occupied by a poroelastic medium, with $\Omega$ and $\Gamma$ an open set and its boundary, respectively.   For the application of boundary conditions, let $\Gamma$ be partitioned as $\Gamma = \overline{\Gamma_u^D \cup \Gamma_u^N} = \overline{\Gamma_f^D \cup \Gamma_f^N}$, with $\Gamma_u^D \cap \Gamma_u^N = \Gamma_f^D \cap \Gamma_f^N = \emptyset$. Here, subscripts $u$ and $f$ denote solid displacement and fluid flow boundary conditions, respectively, while superscripts $D$ and $N$ indicate Dirichlet or Neumann conditions.  Let $\vec{n}_{\Gamma}$ denote the outer normal vector.     Wetting and non-wetting fluid phases are identified by subscript $w$ and $nw$, respectively.   The requisite set of governing equations consists of one conservation law for linear momentum and two conservation laws for the mass balance of each fluid phase.   Let $\vec{u}$, $s_\ell$, and $p_\ell$ denote the solid displacement vector, saturation, and pressure of fluid phase $\ell = \{ w, nw \}$.   Since the two fluids fill the voids, the saturations must satisfy the closure condition $\sum_{\ell = \{w,nw\}} s_\ell = 1$.   In the following, the wetting fluid phase saturation is selected as a primary unknown and will be denoted by lower case $s$ without subscript.   Capillarity effects are not considered, leading to the simplification $p_w = p_{nw} = p$.   This is a frequent assumption in many reservoir engineering applications, but we emphasize that it can have a significant impact on the nature of the governing equations and  solver strategy.  

Given a finite time interval $\mathbb{I} = (0, t_{\max}]$ the strong form of the initial/boundary value problem may be stated as follows:
\begin{subequations}\label{eq:IBVP_global}
\begin{align}
  \intertext{Find $\vec{u} : \overline{\Omega} \times \mathbb{I} \rightarrow \mathbb{R}^3$, $s : \overline{\Omega} \times \mathbb{I} \rightarrow \mathbb{R}$, and $p : \overline{\Omega} \times \mathbb{I} \rightarrow \mathbb{R}$ such that }
  &\nabla \cdot \tensorTwo{\sigma} + \rho\vec{g} = \vec{0} & &\mbox{ on } \Omega \times \mathbb{I} & &\mbox{(quasi-static linear momentum balance),} \label{eq:momentumBalanceS}\\
  %
  &\frac{\partial m_{w}}{\partial t} + \nabla \cdot \vec{w}_{w} = q^I_{w} - q^P_{w} & &\mbox{ on } \Omega \times \mathbb{I} & &\mbox{(mass balance for the wetting fluid phase),} \label{eq:massBalanceW_S}	 \\
  %
  &\frac{\partial m_{nw}}{\partial t} + \nabla \cdot \vec{w}_{nw} = q^I_{nw} - q^P_{nw} & &\mbox{ on } \Omega \times \mathbb{I} & &\mbox{(mass balance for the non-wetting fluid phase),} \label{eq:massBalanceNW_S}	 \\
  %
  \intertext{subject to boundary conditions}
  &\vec{u} = \bar{\vec{u}} & &\mbox{ on } \Gamma_u^D \times \mathbb{I} & &\mbox{(prescribed displacement),} \label{eq:momentumBalanceS_DIR}\\
  &\tensorTwo{\sigma} \cdot \vec{n}_{\Gamma} = \bar{\vec{t}} & &\mbox{ on } \Gamma_u^N \times \mathbb{I} & &\mbox{(prescribed total traction),} \label{eq:momentumBalanceS_NEU}\\
  &p =\bar{p} & &\mbox{ on } \Gamma_f^D \times \mathbb{I} & &\mbox{(prescribed pore pressure),} \label{eq:massBalanceS_DIR}\\    
  &s =\bar{s} & &\mbox{ on } \Gamma_f^D \times \mathbb{I} & &\mbox{(prescribed wetting phase saturation),} \label{eq:massBalanceS_DIR_2}\\        
  &\vec{w}_{w} \cdot \vec{n}_{\Gamma} = \bar{w}_w & &\mbox{ on } \Gamma_f^N \times \mathbb{I} & &\mbox{(prescribed wetting phase mass flux),} \label{eq:massBalanceW_S_NEU}\\           
  &\vec{w}_{nw} \cdot \vec{n}_{\Gamma} = \bar{w}_{nw} & &\mbox{ on } \Gamma_f^N \times \mathbb{I} & &\mbox{(prescribed non-wetting phase mass flux),} \label{eq:massBalanceNW_S_NEU}\\ 
  %
  \intertext{and initial conditions}
  &\vec{u}(\vec{x}, 0) = \vec{u}_0 (\vec{x}) & &\mbox{ } \vec{x} \in \Omega & &\mbox{(initial displacement),} \label{eq:massBalanceDisp_IC}\\
  &s(\vec{x}, 0) = s_0 (\vec{x}) & &\mbox{ } \vec{x} \in \Omega & &\mbox{(initial wetting phase saturation),} \label{eq:massBalanceSatW_IC}\\
  &p(\vec{x}, 0) = p_0 (\vec{x}) & &\mbox{ } \vec{x} \in \Omega & &\mbox{(initial pore pressure).} \label{eq:massBalancePres_IC}
\end{align}
\end{subequations}

\noindent
In Eq. \eqref{eq:momentumBalanceS}, the Cauchy total stress tensor is given by
\begin{equation}
\tensorTwo{\sigma} = \tensorFour{C}_{\text{dr}} : \nabla^s \vec{u} - b p \tensorTwo{1},
\end{equation}
with $\tensorFour{C}_{\text{dr}}$ the fourth-order drained elasticity tensor, $b$ the Biot coefficient, and $\tensorTwo{1}$ the second-order unit tensor; $\rho\vec{g}$ is a body force due to the self-weight of the multiphase mixture, with $\vec{g}$ the gravity vector; $\rho =  [(1-\phi) \rho_s + \phi \rho_w s + \phi \rho_{nw} (1-s) ]$ is the density of the mixture, with $\rho_s$, $\rho_w$, and $\rho_{nw}$ the density of the solid, the wetting, and the non-wetting fluid phase, respectively; $\phi$ is the porosity; and $\nabla^s$  and $\nabla \cdot$ are the symmetric gradient and divergence operator, respectively.

In Eqs. \eqref{eq:massBalanceW_S}, $m_w = (\phi \rho_w s)$, $\vec{w}_w = (\rho_w \vec{v}_w)$, and $q^\alpha_w \ge 0$, $\alpha=\{I,P\}$, denote mass per unit volume, mass flux, and mass injection ($I$) and production ($P$) sources per unit volume for the wetting phase.   Here, $\vec{v}_w$ is the volumetric flux.  It follows the traditional multiphase flow extension of Darcy's law \cite{MusMer36}, 
\begin{equation}
\vec{v}_w = - \lambda_w \tensorTwo{\kappa} \cdot \nabla ( p + \rho_w g z ),
\end{equation}
with $\lambda_w = k_{rw}/\mu_w$ the phase mobility, $\mu_w$ the viscosity, $k_{rw}$ the relative permeability factor,  $\tensorTwo{\kappa}$ the absolute permeability tensor, $g$ the gravitational acceleration, $z$ the elevation above a datum, and $\nabla$ the gradient operator.   In Eqs. \eqref{eq:massBalanceNW_S}, $m_{nw} = (\phi \rho_{nw} (1-s))$, $\vec{w}_{nw} = (\rho_{nw} \vec{v}_{nw})$, $q^I_{nw} \ge 0$, and $q^P_{nw} \ge 0$ are the counterpart of $m_{w}$, $\vec{w}_{w}$, $q^I_{w}$, and $q^P_{w}$ for the non-wetting fluid phase.  

The formulation is completed with constitutive equations and equations of state to specify the following dependencies: $\phi= \phi(\vec{u}, p)$, $\rho_s = \rho_s(p)$, $\rho_\ell = \rho_\ell(p)$, $\mu_\ell = \mu_\ell(p)$, $k_{r\ell} = k_{r\ell} (s)$, $\ell = \{w, nw \}$.   The porosity change from a reference porosity $\phi_0$ is modeled in  incremental form as \cite{Cou04}
\begin{align}
	\Delta \phi &= b \Delta \epsilon_v + \frac{(b - \phi_0)(1-b)}{K\sub{dr}} \Delta p, \label{eq:linearized_porosity}
\end{align}
\noindent
where $\epsilon_v = \text{trace}(\nabla^s \vec{u})$ is the volumetric strain, and $K\sub{dr}$ is the drained (skeleton) bulk modulus.   Specific forms for the remaining relative permeability and fluid compressibility relationships used in our simulations are detailed with the numerical examples in Section \ref{sec:numerical_results}.   The absolute permeability $\tensorTwo{\kappa}$ is here modeled as a constant field and exhibits no stress-dependence, though this additional form of solid-fluid coupling is of frequent interest in applications and often motivates the desire for a fully-coupled solver in the first place.

\section{Discrete Formulation}
\label{sec:num_model}

A variety of choices are available for the spatial discretization of the PDEs governing multiphase poromechanics.   Continuous Finite Element (FE) based formulations---the most popular technique employed in consolidation modeling, e.g. \cite{MurLou94,LewSch98,WhiBor08,Rod_etal16}---have been traditionally used for solving Biot's equations.   However, a discontinuous interpolation for the pressure field may be more appropriate in the presence of highly heterogeneous permeability fields.   In addition, element-wise mass conservation can also be a fundamental requirement, in particular in geoscience applications.   Therefore, formulations that combine FE approaches for mechanics with a Finite Volume (FV) method for fluid flow and transport are often adopted by modelers \cite{SetMou98,KimTchJua11a,Pre14,GarKarTch16}.   Some authors have recently started investigating face-centered and cell-centered finite-volume methods also for the mechanical subproblem \cite{DiP_etal11,Nor14,BofBotDiP16,Nor16,KeiNor17} to account for mechanical effects when flow processes are of predominant importance.   Also, in the context of single-phase poromechanics, locally mass-conservative approaches have been proposed using mixed three-field (displacement-velocity-pressure) or four-field (stress tensor-displacement-velocity-pressure) formulations, respectively \cite[e.g.,][]{Fri00,Lip02,PhiWhe07a,PhiWhe07b,JhaJua07,FerCasGam10,HagOsnLan12b,TurArb14,Yi14,Bot_etal17,BauRadKoc17,Hu_etal17,Rod_etal17,HonKrau18,Bot_etal18}.

We partition the domain using a computational mesh $\mathcal{T}^h$.   The mesh consists of disjoint, nonoverlapping cells such that $\Omega \approx \Omega^h = \bigcup_i K_i$ with $K_i \in \mathcal{T}^h$.   A unique direction is associated with each mesh face $f$ through a unit normal vector $\vec{n}_f$.   For an internal face we let $\vec{n}_f$ point in the direction of increasing global cell indices---i.e. for $f$ shared by cells $K$ and $L$, with $L$ greater than $K$, $\vec{n}_f = \vec{n}_K$ where $\vec{n}_K$ is the outward normal unit vector to cell $K$'s boundary.   For a boundary interface we assume $\vec{n}_f = \vec{n}_{\Gamma}^h$.   In the case of a non-planar face, a plane characterized by a mean normal $\vec{n}_f$ is associated to $f$. 

The numerical solution of the system of governing equations \eqref{eq:IBVP_global} relies on a mixed FE-FV spatial discretization and a fully-implicit time-marching scheme (backward Euler) \cite{Hug00,EymGalHer00,Qua14}.   Without loss of generality, we assume homogeneous prescribed displacement ($\bar{\vec{u}} = \vec{0}$) and wetting and non-wetting fluid phase mass fluxes ($\bar{w}_w = \bar{w}_{nw} = 0$).   Time integration is carried out by partitioning the time interval $\mathbb{I}$ into $n_{\Delta t}$ subintervals $\mathbb{I}_n = (t_{n-1}, t_n]$, $n = 1, \ldots, n_{\Delta t}$, with $\Delta t_n = (t_n - t_{n-1})$ the timestep size.   The mass time derivatives in Eqs. \eqref{eq:massBalanceW_S}-\eqref{eq:massBalanceNW_S} are discretized using discrete mass increments, while the remaining terms are approximated at the end-of-step time $t_n$.   This leads to the following mesh-dependent, fully discrete, weak version  of \eqref{eq:IBVP_global}:
\begin{subequations} \label{eq:weak_form_discr}
\begin{align}
  \intertext{Find $\{ \vec{u}^h_n, s^h_n, p^h_n \} \in \boldsymbol{\mathcal{U}}^h \times \mathcal{S}^h \times \mathcal{P}^h$ such that for $n = \{1, \ldots, n_{n_{\Delta t}} \}$}
  \mathcal{R}_{u} =&\;
	\int_{\Omega^h} \nabla^{s} \tensorTwo{\eta} : \tensorTwo{\sigma}^\prime_n \, \mathrm{d}V -
	\int_{\Omega^h}  \nabla \cdot \tensorTwo{\eta} b p^h_n \, \mathrm{d}V -
	\int_{\Omega^h} \tensorTwo{\eta} \cdot \rho_n \vec{g} \, \mathrm{d}V -
	\int_{\Gamma_u^{N,h}} \tensorTwo{\eta} \cdot \bar{\vec{t}}_n \, \mathrm{d}A =
	0 &&\forall \tensorTwo{\eta} \in \boldsymbol{\mathcal{U}}^h,
	\label{eq:weak_form_mom} \\
	\mathcal{R}_s =&\;
	\int_{\Omega^h} \psi \frac{m_{w,n} - m_{w,(n-1)}}{\Delta t_n}  \, \mathrm{d}V -
  \sum_{f \notin \Gamma_f^{N,h}} \llbracket \psi \rrbracket_f F_{w,n}^f -
  \int_{\Omega^h} \psi (q^I_{w,n} - q^P_{w,n}) \, \mathrm{d}V = 
  0 &&\forall \psi \in \mathcal{S}^h,
  \label{eq:weak_form_massW} \\
	\mathcal{R}_{p} =&\;
	\int_{\Omega^h} \chi \frac{m_{nw,n} - m_{nw,(n-1)}}{\Delta t_n}  \, \mathrm{d}V -
  \sum_{f \notin \Gamma_f^{N,h}} \llbracket \chi \rrbracket_f F_{nw,n}^f -
  \int_{\Omega^h} \chi (q^I_{nw,n} - q^P_{nw,n}) \, \mathrm{d}V = 
  0 &&\forall \chi \in \mathcal{P}^h,
  \label{eq:weak_form_massNW} 
\end{align}
\end{subequations}
\noindent 
using the discrete spaces,
\begin{equation}
	\begin{dcases}
		\boldsymbol{\mathcal{U}}^h := \left\lbrace \tensorTwo{\eta} \left|\right. \tensorTwo{\eta} \in [C^0(\overline{\Omega})]^3, \tensorTwo{\eta} = \vec{0} \text{ on } \Gamma_u^{D,h}, \tensorTwo{\eta}_{\left| K \right.} \in [\mathbb{Q}_1(K)]^3 \; \forall K \in \mathcal{T}^h \right\rbrace, \\
		\mathcal{S}^h := \left\lbrace \psi \left|\right. \psi \in L^2(\Omega), \psi_{\left| K \right.} \in \mathbb{P}_0(K) \; \forall K \in \mathcal{T}^h  \right\rbrace,	\\
		\mathcal{P}^h := \left\lbrace \chi \left|\right. \chi \in L^2(\Omega), \chi_{\left| K \right.} \in \mathbb{P}_0(K) \; \forall K \in \mathcal{T}^h  \right\rbrace,	
	\end{dcases}
	\label{eq:func_discrete_spaces}
\end{equation}
\noindent
with $C^0(\overline{\Omega})$ and $L^2(\Omega)$ the space of continuous and square Lebesgue-integrable functions on $\overline{\Omega}$ and $\Omega$, respectively, $\mathbb{Q}_1(K)$ the space of trilinear polynomials in $K$, and $\mathbb{P}_0(K)$ the piecewise constant space in $K$.   In \eqref{eq:weak_form_massW}-\eqref{eq:weak_form_massNW}, $F_{\ell,n}^f$ denotes a discrete approximation to the inter-cell mass flux for the $\ell$ fluid phase, such that $F_{\ell,n}^f \approx -\int_f \vec{w}_{\ell,n} \cdot \vec{n}_f \; \mbox{d}\Gamma$, with the summation taken over faces in $\mathcal{T}^h$ that do not belong to the (impervious) boundary $\Gamma_f^{N,h}$.   The symbol $\llbracket \cdot \rrbracket_f$ denotes the jump of a quantity across a face $f$ in $\mathcal{T}^h$.   For an internal face, $\llbracket \chi \rrbracket_f = ( {\chi}_{\left| L \right.} - {\chi}_{\left| K \right.} )$, with ${\chi}_{\left| L \right.}$ and ${\chi}_{\left| K \right.}$ the restriction of $\chi$ on cells $K$ and $L$ sharing $f$, respectively.   For a face belonging to the domain boundary, the jump expression reduces to $\llbracket \chi \rrbracket_f = - {\chi}_{\left| K \right.}$.

\noindent 

Discrete approximations for the displacement, pressure, and saturation field are given by
\begin{align}
	\vec{u}^h_n(\vec{x}) &= \sum_{i=1}^{n_u} \tensorTwo{\eta}_i(\vec{x}) u_{i,n}, &
	s^h_n(\vec{x}) &= \sum_{j=1}^{n_s} \psi_j(\vec{x}) s_{j,n}, &
	p^h_n(\vec{x}) &= \sum_{k=1}^{n_p} \chi_k(\vec{x}) p_{k,n},
   \label{eq:uqS_approx}
\end{align}
\noindent
with $\{ \tensorTwo{\eta}_i \}$, $\{ \psi_j \}$, and $\{ \chi_k \}$ bases for $\boldsymbol{\mathcal{U}}^h$, $\mathcal{S}^h$, and $\mathcal{P}^h$, respectively.   The nodal displacement components $\{u_{i,n}\}$, cell-centered wetting phase saturations $\{s_{j,n}\}$, and cell-centered pressures $\{p_{k,n}\}$ are collected in algebraic vectors $\Vec{u}_n \in \mathbb{R}^{n_u}$, $\Vec{s}_n \in \mathbb{R}^{n_s}$, and $\Vec{p}_n \in \mathbb{R}^{n_p}$.   

\begin{rem}
It should be noted that the $\mathbb{Q}_1-\mathbb{P}_0-\mathbb{P}_0$ space adopted here is not automatically inf-sup stable. If the fluid and solid constituents are nearly incompressible, and if the permeability or timestep are very small, the resulting discrete system may allow for non-physical pressure modes, and pressure checkerboarding may be observed.   To address this issue, we include additional pressure jump stabilization terms to the mass balance equations following the macroelement stabilization scheme proposed in \cite{KekSil92}.  Techniques for stabilizing mixed discretizations of multiphase poromechanics are not well studied to date, however, and will be the subject of a separate contribution.  Whenever sufficiently large timesteps or compressible constituents are used, spurious pressure modes are avoided, and such conditions are frequently satisfied in practice.
\end{rem}

A linear two-point flux approximation (TPFA) scheme is adopted for the intercell numerical fluxes $F_{\ell,n}^f$.   For an internal face $f$ shared by cells $K$ and $L$ the fluxes are computed as
\begin{align}
  F_{\ell,n}^f &= \rho^{\texttt{upw}}_{\ell,n} \lambda^{\texttt{upw}}_{\ell,n} \Phi_{\ell,n}^f, &
  \Phi_{\ell,n}^f &= \Upsilon^f \left[ \left( p_{L,n} + \varrho^f_{\ell,n} g z_L \right) - \left( p_{K,n} + \varrho^f_{\ell,n} g z_K \right) \right], &
  \ell \in \{ w, nw \},
  \label{eq:num_flux_approx}
\end{align}

\noindent
where the density $\rho^{\texttt{upw}}_{\ell,n}$ and mobility $\lambda^{\texttt{upw}}_{\ell,n}$ are evaluated using single-point upstream weighting (SPU) according to the sign of the potential $\Phi^f_{\ell,n}$, i.e. the so-called \textit{geometric part} of the $\ell$-fluid phase mass flux~\citep{AziSet79}.   In \eqref{eq:num_flux_approx}, $\Upsilon^f$ is the transmissibility coefficient between cell $K$ and $L$, $z_K$ and $z_L$ are the elevation for cell $K$ and $L$ evaluated at the respective centroids, and $\varrho^f_\ell$ is the mass density averaged at face $f$ as \cite{ZhoTchMal11}

\begin{align}
  \varrho^f_{\ell,n} &=  
  \begin{dcases}  
    (\rho_{\ell,n}^K + \rho_{\ell,n}^L) / 2, & \text{if phase $\ell$ appears in both $K$ and $L$}, \\
    \rho_{\ell,n}^K                    , & \text{if phase $\ell$ appears only in $K$}, \\
    \rho_{\ell,n}^L                    , & \text{if phase $\ell$ appears only in $L$}, \\    
    0                              , & \text{otherwise}.
  \end{dcases}
  \label{eq:avg_density}
\end{align}

\noindent
For a boundary face, $F_{\ell,n}^f$ is still computed using Eqs. \eqref{eq:num_flux_approx}-\eqref{eq:avg_density}, with the neighboring cell index $L$ being replaced by the face index.   See~\ref{app:TPFA_computation} for details concerning the computation of the transmissibility coefficient $\Upsilon^f$.

\begin{rem}
The linear TPFA scheme is the industry standard for reservoir simulation because of its robustness, ease of implementation, and monotonicity-preserving properties.   However, the consistency of the two-point scheme is not guaranteed for arbitrary grid and permeability configurations.  Therefore, the linear TPFA approach may lead to inaccurate results for highly distorted grids or full-permeability tensors. In such cases more sophisticated discretization methods like multipoint and/or nonlinear schemes are needed.   For details on recent developments for heterogeneous anisotropic diffusion problems see~\cite{Dro14,TerMalTch17,Sch_etal18} and references therein.
\end{rem}

The source terms in Eqs. \eqref{eq:weak_form_massW}-\eqref{eq:weak_form_massNW} are used to capture the effect of injection and production wells, modeled as line sources.   These sources are governed by a suitable well model---typically derived by assuming radial flow in the immediate vicinity of the wellbore---relating well control parameters, such as bottomhole pressure, to mass flow rates through the wellbore~\citep{Pea78}.   In this work, we assume each well segment is vertical and connected to the centroid of a cell.   In particular, for a cell $K$ connected to a well the source terms are computed as \citep{CheHuaMa06}
\begin{subequations} \label{eq:peaceman}
\begin{align}
    q_{\ell,n}^P &= - \rho_{\ell,n} \lambda_{\ell,n} \Phi^{K,W}_{\ell,n} \delta \left(\vec{x} - \vec{x}_K \right), \label{eq:peaceman_PROD} \\
    q_{\ell,n}^I &= \rho_{\ell,n} \left(\lambda_{w,n} + \lambda_{nw,n} \right) \Phi^{K,W}_{\ell,n} \delta \left(\vec{x} - \vec{x}_K \right), \label{eq:peaceman_INJ}
\end{align}
\end{subequations}
\noindent 
with
\begin{align}
  \Phi^{K,W}_{\ell,n} &= WI [ (\bar p_{bh,n} + \rho_{\ell,n} g z_{bh}) - (p_{K,n} + \rho_{\ell,n} g z_K ) ],
  \label{eq:peaceman2}
\end{align}
\noindent
where $\bar p_{bh,n}$ is the bottom hole pressure, $\delta \left(\vec{x} - \vec{x}_K \right)$ is the Dirac function with $\vec{x}_K$ the cell centroid position vector, and $WI$ is the well index.   In the present work all wells are assumed to operate under specified bottom-hole pressure conditions for simplicity. Clearly, the condition $\Phi^{K,W}_{\ell,n} < 0$ must hold true to have fluid flow from the reservoir into a production well.   Similarly, injection wells require $\Phi^{K,W}_{\ell,n} > 0$.  Note that injection and production wells differ in the treatment of the phase mobilities.  In Eq.~\eqref{eq:peaceman_INJ}, the sum of phase mobilities is used for injector wells so that saturation conditions are reflected in the wellbore.  This is the standard approach adopted in reservoir simulators when no wellbore crossflow occurs~\cite{Ecl13}.  The well index depends on the well properties, local permeability, and cell dimensions.  For a hexahedral $\tensorTwo{\kappa}$-orthogonal grid---i.e., when the principal axes of the permeability tensor are aligned with the grid---$WI$ is given as
\begin{align}
  WI &= \frac{2 \pi h_z \sqrt{\kappa_x \kappa_y}}{\text{log}(r / r_w) + s_k}, &
  r &= 0.28 \frac{[(\kappa_y/\kappa_x)^{1/2} h_x^2 + (\kappa_x/\kappa_y)^{1/2} h_y^2]^{1/2}}{(\kappa_y/\kappa_x)^{1/4} + (\kappa_x/\kappa_y)^{1/4}},
  \label{eq:well_index}
\end{align}
\noindent
where $\kappa_x$, $\kappa_y$, and $\kappa_z$ are the diagonal  (principal) components of $\tensorTwo{\kappa}$; $h_x$, $h_y$, and $h_z$ are the cell dimensions; $r_w$ is the well radius; and $s_k$ is the well skin factor.    

\section{Newton-Krylov Solver}
\label{sec:NewKry_solver}

The discrete form of the multiphase poromechanical problem \eqref{eq:IBVP_global}  is obtained by introducing expressions~\eqref{eq:uqS_approx}, \eqref{eq:num_flux_approx}, and \eqref{eq:peaceman} into the weak form~\eqref{eq:weak_form_discr}.   This leads to a system of nonlinear equations at time $t_n$,
\begin{align} 
   \begin{cases}
     \Vec{r}_n^{u\;} (\blkVec{x}_n) &= \Vec{0},\\
	  \Vec{r}_n^{s\;}	(\blkVec{x}_n,\blkVec{x}_{n-1}) &= \Vec{0},\\
	  \Vec{r}_n^{p\;} (\blkVec{x}_n,\blkVec{x}_{n-1}) &= \Vec{0},
   \end{cases}
   	\label{eq:NL_res}
\end{align}
\noindent
for the latest solution vector $\blkVec{x}_n=\{ \Vec{u}_n, \Vec{s}_n, \Vec{p}_n\}$. The solution $\blkVec{x}_{n-1}$ is known from the previous timestep.   Starting from an initial guess $\blkVec{x}_n^0$, Newton's method is used to drive the norm of the combined residual vector to below a specified relative tolerance, $\| \blkVec{r}^k \| / \| \blkVec{r}^0 \| < \tau$ .  At each iteration $k$, the Newton update is given by
\begin{subequations}
\begin{align}
\blkMat{A}^k \Delta \blkVec{x} = -\blkVec{r}^k,\\
\blkVec{x}_n^k = \blkVec{x}_n^{k-1} + \Delta \blkVec{x} ,
\end{align}
\end{subequations}
with Jacobian matrix $\blkMat{A}^k = \partial \blkVec{r}^k / \partial \blkVec{x}_n$.  In the current context, the linearization of three coupled governing equations leads to a Jacobian system characterized by an inherent $3 \times 3$ block structure, 
\begin{align}
   \begin{bmatrix}
		\Mat{A}\sub{uu} & \Mat{A}\sub{us} & \Mat{A}\sub{up} \\
		\Mat{A}\sub{su} & \Mat{A}\sub{ss} & \Mat{A}\sub{sp} \\
		\Mat{A}\sub{pu} & \Mat{A}\sub{ps} & \Mat{A}\sub{pp}
	\end{bmatrix}^k
	\begin{bmatrix}
	   \Delta \Vec{u}\\ \Delta \Vec{s} \\ \Delta \Vec{p}
	\end{bmatrix}
	&=
	-
	\begin{bmatrix}
	   \Vec{r}_n^{u} \\ \Vec{r}_n^s \\ \Vec{r}_n^{p}
	\end{bmatrix}^k.
	\label{eq:jac_system}
\end{align}
\noindent 
Detailed expressions for the residual vectors and sub-matrices in $\blkMat{A}$ are given in~\ref{app:FEFV_vec_mat_blocks}.

The linear solution of the Jacobian system~\eqref{eq:jac_system} at each Newton iteration is the most memory-demanding and time-consuming computational kernel in a multiphase poromechanics simulation.  It is also the most difficult component to address in a scalable way on large parallel platforms.  In this work, the system~\eqref{eq:jac_system} is solved iteratively using a nonsymmetric Krylov subspace solver---specifically the generalized minimal residual (GMRES) method~\citep{SaaSch86}.   For large problems, practical convergence of a Krylov method is only possible provided that an effective preconditioner is available.  Using right-preconditioning, we introduce a preconditioning operator $\blkMat{M}^{-1}$ and work with the modified system
\begin{subequations}
\begin{align}
\left( \blkMat{A} \blkMat{M}^{-1} \right)^k {\Delta \blkVec{y}}  = - \blkVec{r}^k, \\ \Delta \blkVec{x} = \blkMat{M}^{-1} {\Delta  \blkVec{y}}.
\end{align} 
\end{subequations}
In the next section, we describe an efficient and scalable two-stage preconditioner that can address both the ill-conditioning and coupled nature of the system at hand.  

\begin{rem}
For simplicity of presentation, we have restricted the discussion to a model system with two immiscible phases.  Nevertheless, the approach devised here is naturally extendible to a larger number of phases and components.  In this case, the \emph{saturation} unknowns $\Vec{s}_n$ should be understood to include all non-pressure, cell-based unknowns, e.g. phase compositions.
\end{rem}

\begin{rem}
Multiphase flow introduces strong nonlinearities, and a na{\"i}ve implementation of Newton's method may fail to converge when large timesteps are taken.  More robust simulators include adaptive timestepping and/or globalization strategies to expand the neighborhood of convergence.  Here, we apply a simple backtracking line search whenever a Newton step fails to satisfy the necessary descent conditions.  We remark that smart strategies to address the intrinsic nonlinearities in these systems remain an active area of research.
\end{rem}

\begin{rem}
The system of residual equations (\ref{eq:NL_res}) stem from the discretization of different governing equations that are not dimensionally consistent with one another.  As a result, it is quite easy for the resulting discrete equations to be poorly scaled with respect to one another, causing a variety of problems at the linear and nonlinear solver level.  A good choice of dimensionally consistent units can typically avoid many of these problems.  As an extra safeguard, we also introduce a simple block-row scaling to ensure that the flow and mechanics equations have roughly equal magnitude before entering the solution routines.
\end{rem}

\begin{rem}
The Jacobian systems~\eqref{eq:jac_system} is solved to a desired linear tolerance.  Thus two tolerances must be selected: the Newton tolerance $\tau$ and the Krylov tolerance $\eta$.  While a constant Krylov tolerance can be employed, this choice will often lead to \emph{oversolving} when Newton iterates are far from the converged solution.  That is, significant work will be invested in solving system~\eqref{eq:jac_system} to a stringent tolerance, even though a much cruder solution is often acceptable for early Newton iterates.  A more efficient strategy is to adapt the linear solver tolerance to reflect the observed nonlinear convergence profile---often referred to as an \emph{inexact} Newton method.  In this way, linear solver work is minimized in early Newton iterates.  A stringent linear solver tolerance is then used in later iterates to maintain the rapid nonlinear convergence behavior of Newton's method in a neighborhood of the solution.  In particular, a good approach is the forcing strategy proposed in \cite{EisWal94}.  Given parameters $\gamma \in (0,1]$, $\omega \in (1,2]$, and $\eta_0 \in [0,1]$, set 
\begin{equation}
\eta_k = \gamma \left( \frac{ \| \Vec{r}^k \| }{\| \Vec{r}^{k-1} \| } \right)^\omega,
\end{equation}
where $\eta^k$ is the relative Krylov tolerance adopted at Newton step $k$.  In practice, an additional safeguard is also included to prevent the tolerance from becoming too small far away from the solution,
\begin{equation}
\eta_k \gets \max [ \eta_k, \gamma \eta_{k-1}^\omega ] \quad \text{whenever} \quad \gamma \eta_{k-1}^\omega > 0.1.
\end{equation}
The specific choices $\gamma=0.9$, $\omega = 2$, and $\eta_0=10^{-1}$ work well in practice.  As the primary focus of the current work is the effectiveness of the linear solver and proposed preconditioning strategy, however, in the numerical examples below we disable this adaptive tolerance and instead set a fixed linear solver tolerance of $\eta=10^{-6}$.    This stringent tolerance leads to oversolving, but it makes the assessment of linear solver performance more straightforward.
\end{rem}

\section{Two-Stage Preconditioner}
\label{sec:precond}

The preconditioner operator $\blkMat{M}^{-1}$ for the Jacobian matrix $\blkMat{A}$ is constructed using two ideas.  The first component is a particular block factorization of the system matrix, which allows us to break the coupled multi-physics problem into simpler sub-problems.  The second component is sparse Schur-complement approximation, allowing us to circumvent certain dense operators which would otherwise appear during the block factorization process.

\subsection{Mechanics / Flow Partitioning}
For clarity of presentation, it is helpful to first re-partition the Jacobian system into a $2\times2$ form in which flow variables (saturations and pressures) are grouped and separated from the mechanics variables (displacements):
\begin{align}
\blkMat{A} = 
\left[
   \begin{array}{c:cc}
		\Mat{A}\sub{uu} & \Mat{A}\sub{us} & \Mat{A}\sub{up} \\ \hdashline
		\Mat{A}\sub{su} & \Mat{A}\sub{ss} & \Mat{A}\sub{sp} \\
		\Mat{A}\sub{pu} & \Mat{A}\sub{ps} & \Mat{A}\sub{pp}
	\end{array} \right]
=
   \begin{bmatrix}
		\Mat{A}\sub{uu} & \Mat{A}\sub{uf}  \\
		\Mat{A}\sub{fu} & \Mat{A}\sub{ff}  
	\end{bmatrix}.
\end{align}
A common approach for constructing preconditioning operators for multiphysics systems is to begin with a block factorization of the matrix \cite{BraPas88,BenGolLie05},
\begin{align}
  \blkMat{A} &=
  \blkMat{L}\blkMat{U} = 
      \begin{bmatrix}
		\Mat{A}\sub{uu} &   \\
		\Mat{A}\sub{fu} & \Mat{S}\sub{ff}  
	\end{bmatrix}
  \begin{bmatrix}
    \Mat{I} & \Mat{A}\sub{uu}^{-1} \Mat{A}\sub{uf} \\
        & \Mat{I} 
  \end{bmatrix} \quad \text{with} \quad \Mat{S}\sub{ff} = \Mat{A}\sub{ff} - \Mat{A}\sub{fu} \Mat{A}\sub{uu}^{-1} \Mat{A}\sub{uf}.
  \label{eq:LU} 
\end{align}
The operator $\Mat{S}\sub{ff}$ is the Schur-complement of $\Mat{A}\sub{uu}$ in $\blkMat{A}$.  In anticipation of developments to come, we remark that this is a ``first-level'' Schur-complement, corresponding to elimination of one field from the original $3\times3$ system. The second-level Schur-complement, stemming from elimination of two fields, will appear later in the discussion.

The guiding idea is to use the inverse of the lower-triangular factor as a template, and seek a preconditioning operator $\blkMat{M}^{-1}\approx \blkMat{L}^{-1}$, where
\begin{align}
  \blkMat{L}^{-1}
  =
  \begin{bmatrix}
    \Mat{A}\sub{uu}^{-1} &  \\
    - \Mat{S}\sub{ff}^{-1} \Mat{A}\sub{fu} \Mat{A}\sub{uu}^{-1} & \Mat{S}\sub{ff}^{-1}
  \end{bmatrix}   . 
\end{align}
This strategy is motivated by the observation that $\blkMat{L}^{-1} \blkMat{A} = \blkMat{U}$.  A property of triangular matrices is that their eigenvalues appear on the diagonal, from which we conclude that the matrix $\blkMat{U}$ has a single eigenvalues $\lambda=1$ with multiplicity $n$.  A direct application of $\blkMat{L}^{-1}$ would therefore provide perfect conditioning.  In a practical setting, however, we can only approximate $\blkMat{L}^{-1}$ and so the construction of a useful scheme requires additional work.

The first challenge arises in the treatment of the Schur-complement operator.  It may be expanded as
\begin{equation}
\Mat{S}\sub{ff} =
\begin{bmatrix}
    \Mat{A}\sub{ss} & \Mat{A}\sub{sp} \\
    \Mat{A}\sub{ps} & \Mat{A}\sub{pp}
  \end{bmatrix}
 -
  \begin{bmatrix}
    \Mat{A}_{su}\Mat{A}_{uu}^{-1}\Mat{A}_{us} & \Mat{A}_{su}\Mat{A}_{uu}^{-1}\Mat{A}_{up} \\
    \Mat{A}_{pu}\Mat{A}_{uu}^{-1}\Mat{A}_{us} & \Mat{A}_{pu}\Mat{A}_{uu}^{-1}\Mat{A}_{up}
  \end{bmatrix}.
\end{equation}
This operator is dense and difficult to work with due to the presence of $\Mat{A}\sub{uu}^{-1}$.  As a first step, we therefore seek to replace the exact $\Mat{S}\sub{ff}$ with a sparse approximation.

A first simplification comes from an examination of the block $\Mat{A}\sub{us}$.  This operator represents the impact of saturation changes on the momentum balance equations.  Specifically, a change in saturation changes the mixture density appearing in equation~(\ref{eq:momentumBalanceS}), which will impact the displacement field.  In general, however, this coupling is quite weak, particularly if the saturation field evolves slowly in time or if phase density differences are small.  For preconditioning purposes, we simply assume $\Mat{A}\sub{us} \approx 0$ and eliminate the dense terms from the first block column of $\Mat{S}\sub{ff}$.

For the remaining dense terms, we rely on a physically-motivated approach that has proven effective for single-phase poromechanics, the \emph{fixed-stress} approximation.  This approach was first introduced as an operator splitting method for sequential solution of poromechanical systems~\cite{SetMou98,KimTchJua11a,KimTchJua13,MikWhe13}.   The sequential method was later reinterpreted in \cite{WhiCasTch16} as a block-preconditioned Richardson iteration involving a particular Schur-complement approximation.  The authors showed that much better performance could be achieved by embedding this same preconditioning strategy within a Krylov iteration.  Here, we consider the multiphase extension of this idea.

Let $\sigma_v = (K_{dr} \epsilon_v - b p )$ denote the mean total stress.   Assuming the mean total stress remains fixed over an increment, i.e. $\Delta \sigma_v = 0$, the volumetric strain can be expressed in terms of the pressure increment as
\begin{align}
   \Delta \epsilon_v^{\textsc{(fs)}} = \frac{b}{K_{dr}} \Delta p.
   \label{eq:FS_hyp}
\end{align}
Substituting expression \eqref{eq:FS_hyp} in the constitutive relationship \eqref{eq:linearized_porosity} during the linearization process removes the explicit dependence of the porosity on the displacement vector, with the resulting discrete mass balance equations uncoupled from the linear momentum balance.    Using this decoupling argument for computing an approximation to the exact Schur complement yields
\begin{align}
  \Mat{S}\sub{ff} \approx \Mat{S}\sub{ff}^{\textsc{(fs)}} =
  \begin{bmatrix}
    \Mat{A}\sub{ss} & \Mat{A}\sub{sp}^{\textsc{(fs)}} \\
    \Mat{A}\sub{ps} & \Mat{A}\sub{pp}^{\textsc{(fs)}} 
  \end{bmatrix}
  , \qquad
  \Mat{A}_{sp}^{(\textsc{fs})} = \Mat{A}_{sp} + \Mat{D}\sub{sp}^{\textsc{(fs)}}, \qquad
  \Mat{A}_{pp}^{\textsc{(fs)}} = \Mat{A}_{pp} + \Mat{D}\sub{pp}^{\textsc{(fs)}},
  \label{eq:Schur_1_FS}
\end{align}
where the entry $(i,j)$ of matrices $\Mat{D}_{sp}^{\textsc{(fs)}}$ and $\Mat{D}_{pp}^{\textsc{(fs)}}$ are given by
\begin{subequations}
\begin{align}
   \left[-\Mat{A}\sub{su}\Mat{A}\sub{uu}^{-1} \Mat{A}\sub{up}\right]_{i,j} \approx
   \left[\Mat{D}_{sp}^{\textsc{(fs)}}\right]_{i,j} =&\;
  \int_{\Omega^h} \psi_i \left( \frac{\partial m_{w}}{\partial \epsilon_v^{\textsc{(fs)}}} \frac{\partial \epsilon_v^{\textsc{(fs)}}}{\partial p} \right)_n \chi_j \, \mathrm{d}V
  &&\forall(i,j) \in \{1, 2, \ldots, n_s \} \times \{1, 2, \ldots, n_p \} ,\label{eq:FS_Dps} \\
   \left[-\Mat{A}\sub{pu}\Mat{A}\sub{uu}^{-1} \Mat{A}\sub{up}\right]_{i,j} \approx
   \left[\Mat{D}_{pp}^{\textsc{(fs)}}\right]_{i,j} =&\;
  \int_{\Omega^h} \chi_i \left( \frac{\partial m_{nw}}{\partial \epsilon_v^{\textsc{(fs)}}} \frac{\partial \epsilon_v^{\textsc{(fs)}}}{\partial p} \right)_n \chi_j \, \mathrm{d}V
  &&\forall(i,j) \in \{1, 2, \ldots, n_p \} \times \{1, 2, \ldots, n_p \} .\label{eq:FS_Dpp}
\end{align}
\label{eq:FS_Dps_Dpp}\null
\end{subequations}
Noting the cell-wise constant interpolation used for both pressure and saturation, the matrices $\Mat{D}_{sp}$ and $\Mat{D}_{pp}$ are diagonal.  Expanding the partial derivatives, the diagonal entries for an element $K_i \in \mathcal{T}^h$ with volume $V_i$ are
\begin{subequations}
\begin{align}
\left[\Mat{D}\sub{sp}^{\textsc{(fs)}}\right]_{i,i} &= V_i  \frac{b^2}{K_{dr}} s_n \, \rho_{w,\,n} ,\\
\left[\Mat{D}\sub{pp}^{\textsc{(fs)}}\right]_{i,i} &= V_i   \frac{b^2}{K_{dr}} (1-s_n) \, \rho_{nw,\,n}.
\end{align}
\end{subequations}
The approximate Schur-complement defined in (\ref{eq:Schur_1_FS}) is therefore sparse, with a storage pattern corresponding to the original finite-volume stencil of the mass balance equations.  The fixed stress assumption only modifies the matrix diagonals associated with mass accumulation terms.  
\begin{rem}
For finite-volume based simulation of multiphase problems that account for a larger number of phases (or compositions), the matrix $\Mat{D}\sub{sp}^{\textsc{(fs)}}$ will have a block-diagonal structure with blocks having dimensions $n_s^K \times 1$, where $n_s^K$ denotes the total number of non-pressure degrees-of-freedom per element.
\end{rem}
\begin{rem}
The relationship between volumetric strain and fluid pressure defined by equation~(\ref{eq:FS_hyp})  is an approximation based on the assumption that the loading path is primarily volumetric. The actual relationship will depend both on material properties \emph{and} on the specific loading conditions for a given problem.  While the modulus $K_{dr}$ will provide good convergence in all cases, it may not be the optimal value for any given problem.  For example, under one-dimensional loading, a uniaxial modulus will more accurately approximate the coupling condition.  See, for example, recent work in \cite{MikWhe13,BauRadKoc17} on using different fixed stress assumptions to improve convergence rates.
\end{rem}
\begin{rem}
An algebraic generalization of the fixed-stress approximation was proposed in \cite{Kle_etal16} based on a row-sum lumping (RSL) strategy.   Using a probing vector $\Vec{e} = \{1, \ldots, 1\}^T \in \mathbb{R}^{n_p}$, the entries of the diagonal matrices defined in \eqref{eq:FS_Dps_Dpp} can be computed as
\begin{align}
   \Mat{D}\sub{sp}^{\textsc{(rsl)}} = - \texttt{diag} ( \Mat{A}\sub{su} \Mat{A}\sub{uu}^{-1} \Mat{A}\sub{up}\Vec{e} ), \qquad
   \Mat{D}\sub{pp}^{\textsc{(rsl)}} = - \texttt{diag} ( \Mat{A}_{pu} \Mat{A}_{uu}^{-1} \Mat{A}_{up} \Vec{e} ).
  \label{eq:FS_RSL}
\end{align}
Here, $\texttt{diag}(\cdot)$ indicates  construction of a diagonal matrix using an input vector. Computing matrices $\Mat{D}\sub{sp}^{\textsc{(rsl)}}$ and $\Mat{D}\sub{pp}^{\textsc{(rsl)}}$ in this manner requires two solves with $\Mat{A}\sub{uu}$ up to a desired accuracy, and four matrix-vector multiplications.   If a good preconditioner for $\Mat{A}\sub{uu}$ is available, it can be used to approximate the action of $\Mat{A}\sub{uu}^{-1}$ in the construction.  Note that it is not necessary to update these vectors at each nonlinear iteration if the relevant Jacobian blocks do not change significantly.   An alternative algebraic construction for the diagonal approximation $\Mat{D}\sub{pp}^{\textsc{(fs)}}$ was presented and analyzed in \cite{CasWhiFer16} for single-phase poromechanics using an element-by-element approach.
\end{rem}

Returning now to the template $\blkMat{M}^{-1}\approx \blkMat{L}^{-1}$, the top level preconditioning operator is defined as
\begin{align}
  \blkMat{M}^{-1}
  =
  \begin{bmatrix}
    \Mat{M}\sub{uu}^{-1} &  \\
    - \Mat{M}\sub{ff}^{-1} \Mat{A}\sub{fu} \Mat{M}\sub{uu}^{-1} & \Mat{M}\sub{ff}^{-1}
  \end{bmatrix}
  =
  \begin{bmatrix}
     \Mat{I} & \\
                 & \Mat{M}\sub{ff}^{-1} 
  \end{bmatrix}
  \begin{bmatrix}
     \Mat{I} & \\
     -\Mat{A}\sub{fu} & \Mat{I}
  \end{bmatrix}
  \begin{bmatrix}
     \Mat{M}\sub{uu}^{-1} & \\  
                             & \Mat{I}
  \end{bmatrix}
  ,  
  \label{eq:blk_tr_prec}    
\end{align}
where two sub-preconditioners have been identified, 
\begin{subequations}
\begin{align}
\Mat{M}\sub{uu}^{-1} &\approx \Mat{A}\sub{uu}^{-1}, \\
\Mat{M}\sub{ff}^{-1} &\approx \left(\Mat{S}\sub{ff}^\textsc{(fs)}\right)^{-1},
\end{align}
\end{subequations}
corresponding to a mechanics preconditioner acting on the displacement components, and a flow preconditioner acting on the combined saturation and pressure components.  We will describe particular forms for them in subsequent sections.   In general, though, we remark that the sparse Schur complement system can be viewed as the discretization of a standard multiphase flow problem \emph{without} mechanics, where the fixed stress assumption provides a particular form for the rock-compressibility term.  Therefore, much of the existing knowledge on preconditioning multiphase flow problems may be brought to bear in designing $\Mat{M}\sub{ff}^{-1}$.  The matrix $\Mat{A}\sub{uu}$ is an elasticity operator, and many existing tools also exist for constructing $\Mat{M}\sub{uu}^{-1}$.  Therefore, the block factorization strategy provides a convenient and modular framework for assembling a coupled multiphase poromechanics solver from highly-tuned component pieces.

\subsection{Mechanics Preconditioner}

The elasticity block $\Mat{A}\sub{uu}$ is an elliptic operator that is well suited to multigrid techniques \cite{Baker:2010}.   Here, we take a simple approach known as the \emph{separate displacement component} approximation, introduced in \cite{AxeGus78}.  If the displacement degrees of freedom are ordered by coordinate direction, then the resulting matrix has its own $3\times3$ structure,
\begin{equation}
\Mat{A}\sub{uu} = 
\left[
   \begin{array}{ccc}
		\Mat{A}\sub{xx} & \Mat{A}\sub{xy} & \Mat{A}\sub{xz} \\ 
		\Mat{A}\sub{yx} & \Mat{A}\sub{yy} & \Mat{A}\sub{yz} \\
		\Mat{A}\sub{zx} & \Mat{A}\sub{zy} & \Mat{A}\sub{zz}
	\end{array} \right].
\end{equation}
Rather than building a preconditioner using the complete matrix, we consider the  sparser approximation
\begin{equation}
\Mat{A}^\textsc{(sdc)}\sub{uu} = 
\left[
   \begin{array}{ccc}
		\Mat{A}\sub{xx} &  &  \\ 
		& \Mat{A}\sub{yy} &  \\
		 & & \Mat{A}\sub{zz}
	\end{array} \right],
\end{equation}
in which the off-diagonal blocks have been eliminated. Using Korn's inequality, one can show that $\Mat{A}^\text{(SDC)}\sub{uu}$ is spectrally equivalent to $\Mat{A}\sub{uu}$ \cite{Bla94,GusLin98}.  Note that this approximation breaks down in the incompressible elasticity limit, Poisson ratio $\nu \to 0.5$, though this is not a concern for most geologic materials.  We then construct an algebraic multigrid preconditioner from the block-diagonal approximation, 
\begin{equation}
\Mat{M}\sub{uu}^{-1} = \texttt{amg}\left(\Mat{A}^\textsc{(sdc)}\sub{uu} \right),
\end{equation}
which can be thought of as three scalar AMG preconditioners, one for each coordinate direction.  More specific details regarding the AMG implementation are summarized with the numerical examples below.

The mechanics preconditioner is the most expensive operator to construct during the solver setup phase.  For an elastic model, however, the matrix $\Mat{A}\sub{uu}$ is constant during the entire simulation.  This AMG preconditioner therefore only needs to be initialized once, and its setup cost is quickly amortized over many solver calls.   This is not the case for the flow preconditioner, which is typically re-initialized within every Newton iteration.

\subsection{Flow Preconditioner}
\label{sec:Schur_CPR}

The preconditioning operator for the approximate Schur complement $\Mat{S}^{\textsc{(fs)}}\sub{ff}$ is built following the \textit{Constrained Pressure Residual} (CPR) approach \cite{Wal83,WalKenLit85,Cao_etal05}.  This is a multistage strategy that combines preconditioners with different features in a multiplicative fashion.  In general, given a matrix $\Mat{B}$ in $\mathbb{R}^{n \times n}$, a multiplicative multistage preconditioner $\Mat{M}_{B}^{-1}$ for $\Mat{B}$ can be written as~\cite{Cao_etal05}

\begin{align}
  \Mat{M}_B^{-1} &= \Mat{M}_1^{-1} + \sum_{i=2}^{n_{\text{st}}} \Mat{M}_i^{-1} \prod_{j=1}^{i-1} \left(\Mat{I} - \Mat{B} \Mat{M}_j^{-1} \right),
  \label{eq:multi_prec}
\end{align}
with $n_{\text{st}}$ the number of stages, $\Mat{M}_i^{-1}$  the $i$th preconditioner for $\Mat{B}$, and $\Mat{I}$ the identity matrix in $\mathbb{R}^{n \times n}$.   

The design of the CPR multistage approach---the industry standard for fully-implicit simulation of multiphase flow---is guided by the underlying nature of the governing equations.   Fields that exhibit near-elliptic behavior, such as pressure, are characterized by long-range (global) coupling.  They require multilevel preconditioners to reduce error modes that are smooth across the computational domain.   Conversely, fields with near-hyperbolic character, such as saturation, display short-range (local) coupling.  Cheaper strategies can be employed to remove the associated error modes.   Standard CPR is implemented as a two-stage method that combines a global (${\Mat{M}}^{-1}_G$) preconditioner with a local preconditioner (${\Mat{M}}^{-1}_L$).    The overall preconditionioning operator $\Mat{M}\sub{ff}^{-1}$ is formally 
\begin{align}
  \Mat{M}\sub{ff}^{-1} = \left[ {\Mat{M}}_G^{-1} + {\Mat{M}}_L^{-1} \left(\Mat{I} - \Mat{S}\sub{ff}^\textsc{(fs)} {\Mat{M}}_G^{-1} \right) \right].
  \label{eq:mult_prec_op}
\end{align}

\noindent
In practice, $\Mat{M}\sub{ff}^{-1}$ is applied through a sequence of vector-multiply operations.


The derivation of  ${\Mat{M}}_G^{-1}$ also relies on a block LU-factorization, this time of the approximate Schur-complement $\Mat{S}\sub{ff}^{\textsc{(fs)}}$, 

\begin{align}
   \Mat{S}\sub{ff}^{\textsc{(fs)}} &= 
   \begin{bmatrix}
      \Mat{A}\sub{ss} & \\
      \Mat{A}\sub{ps} & \Mat{S}\sub{pp}
   \end{bmatrix}
   \begin{bmatrix}
      \Mat{I} & \Mat{A}\sub{ss}^{-1} \Mat{A}\sub{sp}^{\textsc{(fs)}} \\
      & \Mat{I}
   \end{bmatrix}
   \quad \text{with} \quad
   \Mat{S}\sub{pp} = \Mat{A}\sub{pp}^{\textsc{(fs)}} - \Mat{A}_{ps} \Mat{A}_{ss}^{-1} \Mat{A}_{sp}^{\textsc{(fs)}}.
   \label{eq:S1_LDU}
\end{align}
Note that $\Mat{S}\sub{pp}$ is the ``second-level'' Schur complement, obtained after recursive elimination of both displacements and saturations from the original $3\times3$ system $\blkMat{A}$.   We again seek to sparsify the Schur-complement, this time by introducing diagonal approximations for $\Mat{A}_{ss}$ and $\Mat{A}_{ps}$.  Two common strategies are \emph{True-IMPES reduction} and \emph{Quasi-IMPES reduction}.

\textit{True-IMPES reduction} seeks to mimic the \textit{IMplicit-Pressure Explicit-Saturation} (IMPES) time discretization scheme \cite{AziSet79} on the linear level.   A column-sum lumping operation is used to lump off-diagonal entries---containing derivatives of the flux terms with respect to the saturation degrees of freedom---into the diagonal.  This can be computed using the matrix transpose as
\begin{subequations}
\begin{align}
   \Mat{A}_{ss} \approx \Mat{D}_{ss}^{\textsc{(cpr)}} &= \texttt{diag} ( \Mat{A}_{ss}^T \Vec{e} ), \\
   \Mat{A}_{ps} \approx \Mat{D}_{ps}^{\textsc{(cpr)}} &= \texttt{diag} ( \Mat{A}_{ps}^T \Vec{e} ), 
  \label{eq:CPR_CSL}
\end{align}
\end{subequations}
with $\Vec{e} = \{1, \ldots, 1 \} \in \mathbb{R}^{n_s}$.   A simpler alternative is the \emph{Quasi-IMPES reduction}, which extracts the diagonal entries from the respective matrices, without lumping in off-diagonal terms,
\begin{subequations}
\begin{align}
   \Mat{A}_{ss} \approx \Mat{D}_{ss}^{\textsc{(cpr)}} &= \texttt{diagm} ( \Mat{A}_{ss}), \\
   \Mat{A}_{ps} \approx \Mat{D}_{ps}^{\textsc{(cpr)}} &= \texttt{diagm} ( \Mat{A}_{ps}).
\end{align}
\end{subequations}
To clarify notation, $\texttt{diag}(\cdot)$ indicates  construction of a diagonal matrix using an input \emph{vector}, while $\texttt{diagm}(\cdot)$ indicates construction of a diagonal matrix using an input \emph{matrix}.  In practice, the diagonal matrices are simply stored as a vector of their diagonal entries. Note that in a more general  formulation where multiple saturation or compositional variables are present, the CPR reduction leads to block-diagonal matrices with small, dense blocks corresponding to degrees of freedom coupled in each cell.  See \cite{Cao_etal05} for details.

With these diagonal approximations, the second-level Schur complement may be approximated as
\begin{equation}
\Mat{S}_{pp} \approx \Mat{S}_{pp}^\textsc{(cpr-fs)} =  \Mat{A}\sub{pp}^{\textsc{(fs)}} - \Mat{D}_{ps}^\textsc{(cpr)} \left( \Mat{D}_{ss}^{\textsc{(cpr)}}\right)^{-1} \Mat{A}_{sp}^{\textsc{(fs)}}.
\end{equation}

\noindent
Note that this strategy again preserves the matrix sparsity pattern stemming from the original finite-volume stencil.  The result is a near-elliptic pressure system that is well suited to AMG preconditioning.  As before, we then use the inverse of the lower triangular factor appearing in \eqref{eq:S1_LDU} as a preconditioner template.  This leads to the concrete definition of the global preconditioner
\begin{align}
   \Mat{M}_{G}^{-1} & 
  =  
  \begin{bmatrix}
    \Mat{I} & \\
        & \Mat{M}_{pp}^{-1} 
  \end{bmatrix}  
  \begin{bmatrix}
    \Mat{I}              &  \\
    -\Mat{A}_{ps} & \Mat{I}
  \end{bmatrix}  
  \begin{bmatrix}
    \left(\Mat{D}_{ss}^{\textsc{(cpr)}}\right)^{-1} &  \\
                & \Mat{I}
  \end{bmatrix} ,
  \quad \text{with} \quad \Mat{M}^{-1}_{pp} = \texttt{amg}\left( \Mat{S}_{pp}^\textsc{(cpr-fs)} \right).
   \label{eq:S1_prec}
\end{align}
In equation~(\ref{eq:S1_prec}), we note that the matrix $A\sub{ps}$ could also be replaced here with its diagonal approximation.  We have generally found better performance using the complete matrix, however, despite the larger cost for a matrix-vector multiply operation.

The global preconditioner $\Mat{M}_G^{-1}$ corrects low-frequency errors associated with the pressure unknowns, but it only performs a simple (Jacobi) correction to the saturation unknowns using $\Mat{D}_{ss}^{\textsc{(cpr)}}$.   The second-stage preconditioner $\Mat{M}_L^{-1}$ complements this first stage by smoothing out the remaining error modes.  These residual errors are primarily associated with steep saturation gradients, and are highly localized.  As a result, single-level and processor-local preconditioning strategies can be quite effective.  

To construct the second stage, the matrix $\Mat{S}\sub{ff}^\textsc{(fs)}$ is viewed as a monolithic system. Let $\Mat{P}$ be a permutation matrix that reorders the partitioned vector of saturation and pressure unknowns into an interleaved ordering---that is, degrees of freedom are re-ordered as $\{\{s_1,p_1\},\{s_2,p_2\},...,\{s_n,p_n\}\}$.  This permutation creates a sparsity pattern in which dense $2\times2$ blocks appear for each grid cell, with these cell blocks connected by the TPFA finite volume stencil.  This hierarchical structure also appears in the more general case when multiple saturation or compositional variables exist.  The size of the dense blocks corresponds to the total number of unknowns per cell. Block versions of relaxation or incomplete factorization preconditioners are therefore appealing, as dense multiplication and inversion operations can be applied to the small blocks. 

Common strategies for $\Mat{M}_L^{-1}$ are relaxation preconditioners (e.g. Jacobi, Gauss-Seidel) or incomplete factorizations (e.g. ILU(k), ILUT).  In this work, we explore two options.  Option 1 uses several sweeps of hybrid block Gauss-Seidel (\texttt{hbgs}), i.e.
\begin{equation}
\Mat{M}_L^{-1} = \texttt{hbgs}\left( \Mat{P} \Mat{S}\sub{ff}^\textsc{(fs)} \Mat{P}^T \right),
\end{equation}
where the input matrix is the interleaved version of the partitioned flow system.  Option 2 uses one sweep of processor-local, pointwise ILU(k) (\texttt{iluk}) on the interleaved system.   

\subsection{Final Algorithm}

\begin{algorithm}[t]
	\caption{Application of the preconditioner $\blkMat{M}^{-1}$}
	\label{alg:apply_M_inv}
	\begin{algorithmic}[1]
		\Require $\Vec{v} = \begin{bmatrix} \Vec{v}\sub{u} ; \Vec{v}\sub{s} ; \Vec{v}\sub{p} \end{bmatrix} $
		\Ensure $\Vec{z} = \blkMat{M}^{-1} \Vec{v}$
		\State $\Vec{z}\sub{u} = \Mat{M}\sub{uu}^{-1} \Vec{v}\sub{u}$
		\State $\begin{bmatrix} \Vec{y}\sub{s} \\  \Vec{y}\sub{p} \end{bmatrix} =
		\begin{bmatrix} \Vec{v}\sub{s} \\  \Vec{v}\sub{p} \end{bmatrix}
		- \begin{bmatrix} \Mat{A}\sub{su} \\  \Mat{A}\sub{pu} \end{bmatrix} z_u$	
		\State $\Vec{z}\sub{s} = \left(\Mat{D}\sub{ss}^{\textsc{(cpr)}} \right)^{-1} \Vec{y}_s$
		\State $\Vec{w}\sub{p} = \Vec{y}\sub{p} - \Mat{A}\sub{ps}  \Vec{z}\sub{s}$
		\State $\Vec{z}\sub{p} = \Mat{M}\sub{pp}^{-1} \Vec{w}_p$
		\State $\begin{bmatrix} \Vec{w}\sub{s} \\  \Vec{w}\sub{p} \end{bmatrix} =
		\begin{bmatrix} \Vec{y}\sub{s} \\  \Vec{y}\sub{p} \end{bmatrix}
		- \begin{bmatrix} \Mat{A}\sub{ss} & \Mat{A}\sub{sp}^\textsc{(fs)} \\  \Mat{A}\sub{ps} & \Mat{A}\sub{pp}^\textsc{(fs)}\end{bmatrix} 
		\begin{bmatrix} \Vec{z}\sub{s} \\  \Vec{z}\sub{p} \end{bmatrix} $
		\State$  \begin{bmatrix} \Delta \Vec{z}\sub{s} \\ \Delta \Vec{z}\sub{p} \end{bmatrix} = \Mat{P}^T  \Mat{M}_L^{-1} \Mat{P} \begin{bmatrix}  \Vec{w}\sub{s} \\  \Vec{w}\sub{p} \end{bmatrix}$
		\State$  \begin{bmatrix} \Vec{z}\sub{s} \\  \Vec{z}\sub{p} \end{bmatrix} \gets  \begin{bmatrix} \Vec{z}\sub{s} \\  \Vec{z}\sub{p} \end{bmatrix}  +  \begin{bmatrix} \Delta \Vec{z}\sub{s} \\  \Delta \Vec{z}\sub{p} \end{bmatrix} $
	\end{algorithmic}
\end{algorithm}

Using the developments above, the concrete implementation of the preconditioning strategy can be summarized as follows. In the simulator, linear solver calls are broken into two phases: a setup phase and an iteration phase.  During the setup phase, auxiliary vectors and matrices associated with the fixed stress and CPR reductions are assembled.  Then sub-preconditioners $\Mat{M}\sub{uu}^{-1}$, $\Mat{M}\sub{pp}^{-1}$, and $\Mat{M}\sub{L}^{-1}$ can be constructed.  The elasticity preconditioner $\Mat{M}\sub{uu}^{-1}$ is the most expensive to build, but it only needs to be initialized once during the entire simulation.  Once the necessary operators are assembled, a non-symmetric Krylov solver (e.g. GMRES) may be called with a preconditioning routine implementing the strategy proposed here. As an aid to the reader, Algorithm~\ref{alg:apply_M_inv} summarizes the concrete steps necessary to apply the preconditioner $\blkMat{M}^{-1}$ to an input vector.  We observe that steps 1-5 implement a block-lower-triangular preconditioner involving two AMG sub-preconditioners and one Jacobi (diagonal) preconditioner.  Steps 6-7 then implement a second stage correction to the saturation and pressure fields using the cheap local smoother.

\begin{rem}
While this work focuses on ``monolithic'' solver strategy using a preconditioned Krylov iteration, we remark that a ``sequential'' solver strategy can be obtained easily by simply embedding the proposed preconditioner within a Richardson iteration.  Step 1 of the preconditioner implements a correction to the displacement field, while steps 2-8 implement a subsequent update to the flow variables.  One would then iterate back and forth between mechanics and flow in a sequential fashion within the fixed point iteration.  See \cite{WhiCasTch16} for further details.  In general, however, the Krylov approach will provide superior convergence properties.
\end{rem}

\section{Numerical Examples}
\label{sec:numerical_results}

  \begin{table}[p]
\caption{Parameter values used for the numerical examples.  Note that the units presented below are chosen for clarity but do not form a self-consistent unit system.  They are converted to consistent values prior to use within the simulator.}

\label{table:parameters}
\centering
\small
\begin{tabular}{llrr}
\toprule
Parameter & Units & Staircase Example & SPE10-based Example \\
\midrule 
\it Porosity: \\
$\qquad$Reservoir & -- & 0.20 &  {\footnotesize SPE10} data \cite{ChrBlu01} \\
$\qquad$Seal units & -- & 0.05    & 0.01\\[5pt]
\it Permeability: \\
$\qquad$Reservoir & mD & 1000 &  {\footnotesize SPE10} data \cite{ChrBlu01} \\
$\qquad$Seal units & mD & 1    &  0.01\\[5pt]
\it Relative perm: \\
$\qquad$Residual water sat. & -- & 0.2 & 0.2 \\
$\qquad$Residual oil sat. & -- & 0.2 & 0.2 \\[5pt]

\it Water:\\
$\qquad$Reference density &kg/m$^3$& 1035 & 1035 \\
$\qquad$Compressibility &MPa$^{-1}$& $4.34\cdot 10^{-4}$ & $4.34\cdot 10^{-4}$ \\
$\qquad$Viscosity & cP & 0.3  & 0.3 \\[5pt]


\it Oil:\\
$\qquad$Reference density &kg/m$^3$& 863 & 863 \\
$\qquad$Compressibility &MPa$^{-1}$& $1.98\cdot 10^{-4}$ & $1.98\cdot 10^{-4}$ \\
$\qquad$Viscosity & cP & 3.0  & 3.0 \\[5pt]

\it Rock: \\
$\qquad$Young's modulus &MPa& 5000 & 5000 \\
$\qquad$Biot coefficient & -- & 1 & 1  \\
$\qquad$Grain density & kg/m$^3$ & 2650 & 2650\\[5pt]

\it Well control:\\
$\qquad$Injection $\Delta$BHP &MPa& $5$ & $5$ \\
$\qquad$Production $\Delta$BHP &MPa& $-5$ & $-5$ \\
$\qquad$Ramp time & day & 1 & 1\\
$\qquad$Well radius & m & 0.1524 & 0.1524\\
$\qquad$Skin factor & -- & 0 & 0\\[5pt]
\it Time-stepping: \\
$\qquad$Initial $\Delta t$ & day & 0.1 & 0.1 \\
$\qquad$Maximum $\Delta t$ & day & 1 & 5\\
$\qquad$End time & day & 100 & 100\\[5pt]

\it Solver tolerances: \\
$\qquad$Newton & -- & 10$^{-5}$ & 10$^{-5}$ \\
$\qquad$Krylov  & -- & 10$^{-6}$ & 10$^{-6}$ \\
\bottomrule
\end{tabular}
\end{table}

 We present two numerical examples to test the performance of the proposed preconditioner: (1) a weak scaling study using a simple, synthetic configuration; and (2) a strong scaling study using a more realistic, highly-heterogeneous reservoir. The code used for this study, \emph{Geocentric}, relies heavily on the \emph{deal.ii} Finite Element Library \citep{BanHarKan07} for discretization functionality.  Algebraic multigrid, relaxation, and incomplete factorization sub-preconditioners are provided by \emph{Hypre} \cite{FalYan02}.   All numerical results were obtained on a cluster of Intel Xeon E5-2695 processors.  Each computational node contains two 18-core processors sharing 128 GiB of memory, with Intel Omni-Path interconnects between nodes.  We use pure MPI-based parallelism.
 

In general, the performance of the proposed algorithm is heavily dependent on the specific design of the sub-preconditioners.  The numerical results below are based on the following choices: For the separate-displacement-component elasticity operator, we use an AMG preconditioner with one level of aggressive coarsening, one sweep of hybrid forward $\ell_1$-Gauss-Seidel \cite{Baker:2011} for the down cycle and one sweep of hybrid backward $\ell_1$-Gauss-Seidel for the up cycle. The coarsest grid is solved directly with Gaussian elimination. For the first-stage flow system, we use quasi-IMPES reduction.  The pressure system uses a standard scalar AMG approach with a Hybrid Modified Independent Set (HMIS) coarsening strategy \cite{DeSterck06}. Similar to the elasticity block, we use one sweep of hybrid forward $\ell_1$-Gauss-Seidel for the down cycle, one sweep of hybrid backward $\ell_1$-Gauss-Seidel for the up cycle, and a direct solve with Gaussian elimination for the coarsest grid. For the second stage correction, we explore two options.  For the first benchmark in \S\ref{ex:staircase} we use three sweeps of hybrid block Gauss-Seidel.  For the second benchmark in \S\ref{ex:spe10} we use one sweep of processor-local, pointwise ILU(0).

Details of the material and parameter values used in both examples are summarized in Table~\ref{table:parameters}. In both examples, fluid phases follow the density model
\begin{equation}
\rho_l = \rho^o_l \left[ 1 + c_l \left(p-p^o\right) \right],
\end{equation}
where $c_l$ is the fluid compressibility and $\rho_l^o$ is the phase density at a reference pressure  $p^o$.  This reference pressure is taken as initial formation pressure.  For the rock, solid grains are intrinsically incompressible with Biot coefficient $b=1$, maximizing the coupling strength between the mechanical and flow problems.  The relative permeability relationships are
\begin{equation}
k_{rw} = \hat{s}^2, \qquad k_{rnw} = \left(1-\hat{s}\right)^2,\qquad \hat{s} = \frac{s-s_{wr}}{1-s_{wr}-s_{nwr}},
\end{equation}
where $s_{wr}$ and $s_{nwr}$ are the wetting and non-wetting residual saturations.  At $t=0$, the formations are assumed to be in hydrostatic and lithostatic equilibrium, with a uniform saturation $s(\vec{x},0) = s_{wr}$.  In both examples, target well pressures are not applied instantaneously, but rather ramped up linearly over a fixed time period (1 day) to relax the stiffness of the nonlinear problem.  We use a telescoping time-stepping scheme, where simulation timesteps grow exponentially from an initial value to a maximum value, after which they remain constant for the remainder of the simulation.

\subsection{Staircase Benchmark}
\label{ex:staircase}

\begin{figure*}
    \centering
    \hfill
    \begin{subfigure}[t]{0.48\textwidth}
        \includegraphics[width=\textwidth]{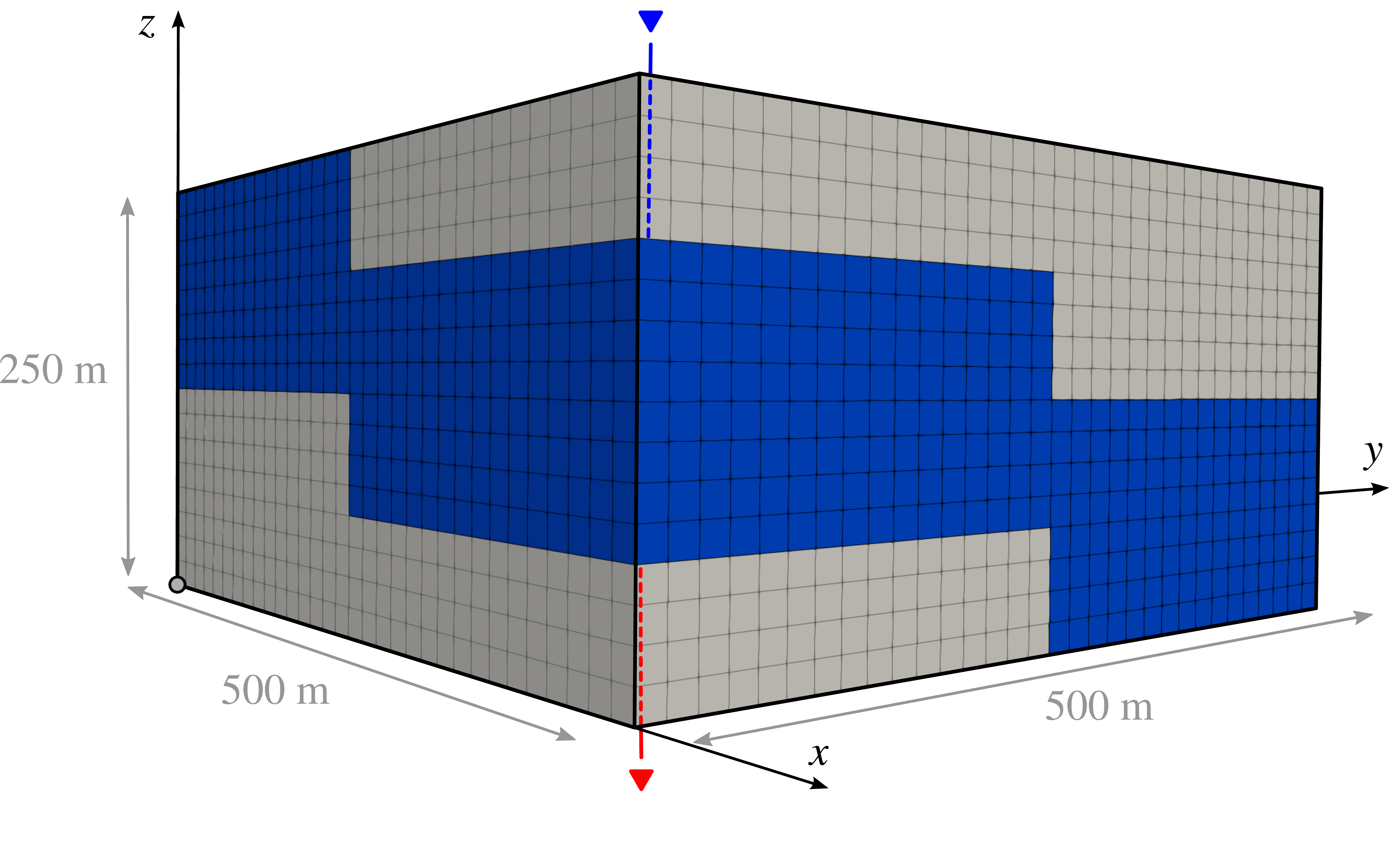}
        \caption{}
        \label{fig:staircase_domain}
    \end{subfigure}
    \hfill
    \begin{subfigure}[t]{0.48\textwidth}
        \includegraphics[width=\textwidth]{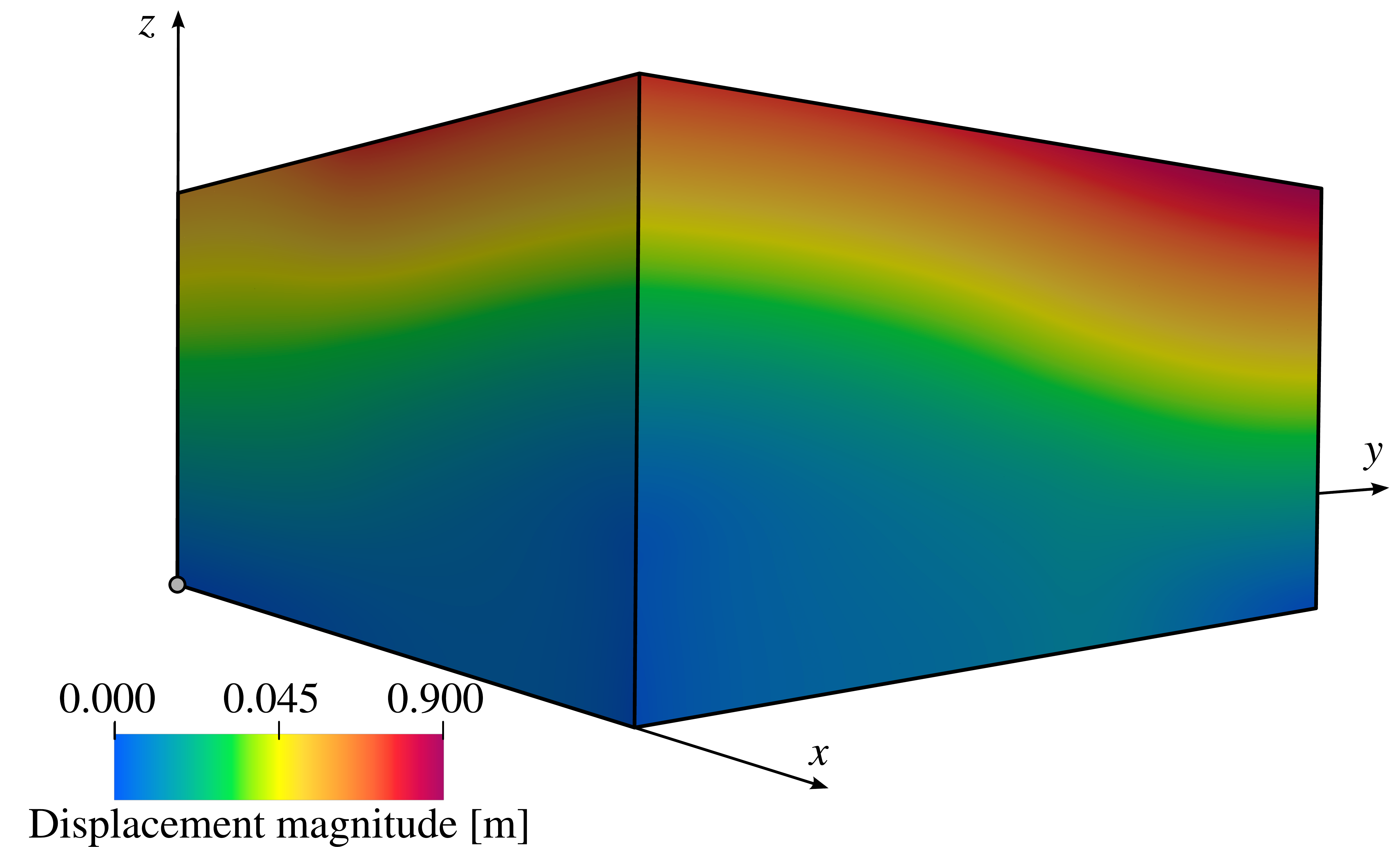}
        \caption{}
        \label{fig:staircase_disp}
    \end{subfigure}
    \hfill\null
    
    \hfill
    \begin{subfigure}[t]{0.48\textwidth}
        \includegraphics[width=\textwidth]{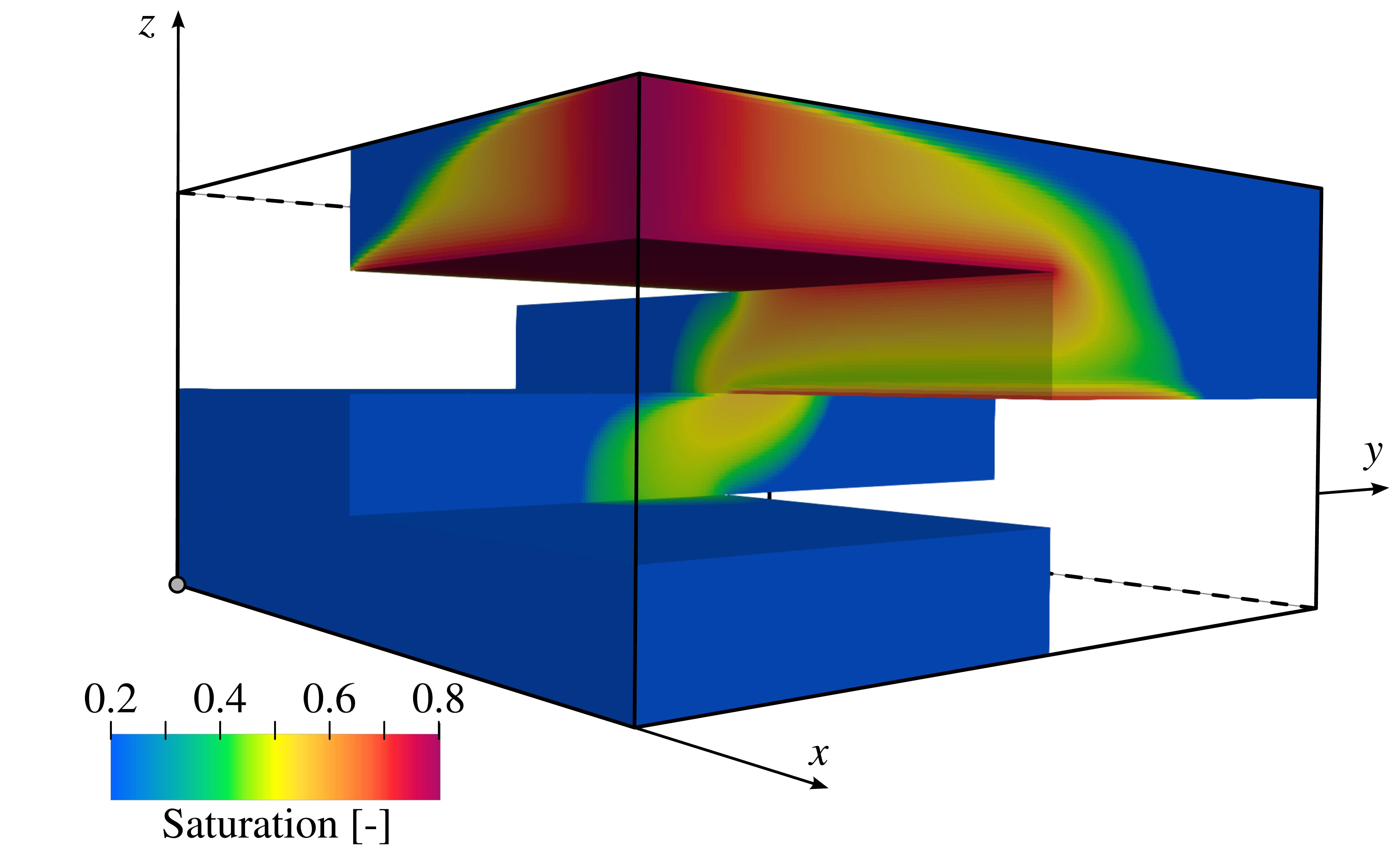}
        \caption{}
        \label{fig:staircase_sat}
    \end{subfigure}   
    \hfill
    \begin{subfigure}[t]{0.48\textwidth}
        \includegraphics[width=\textwidth]{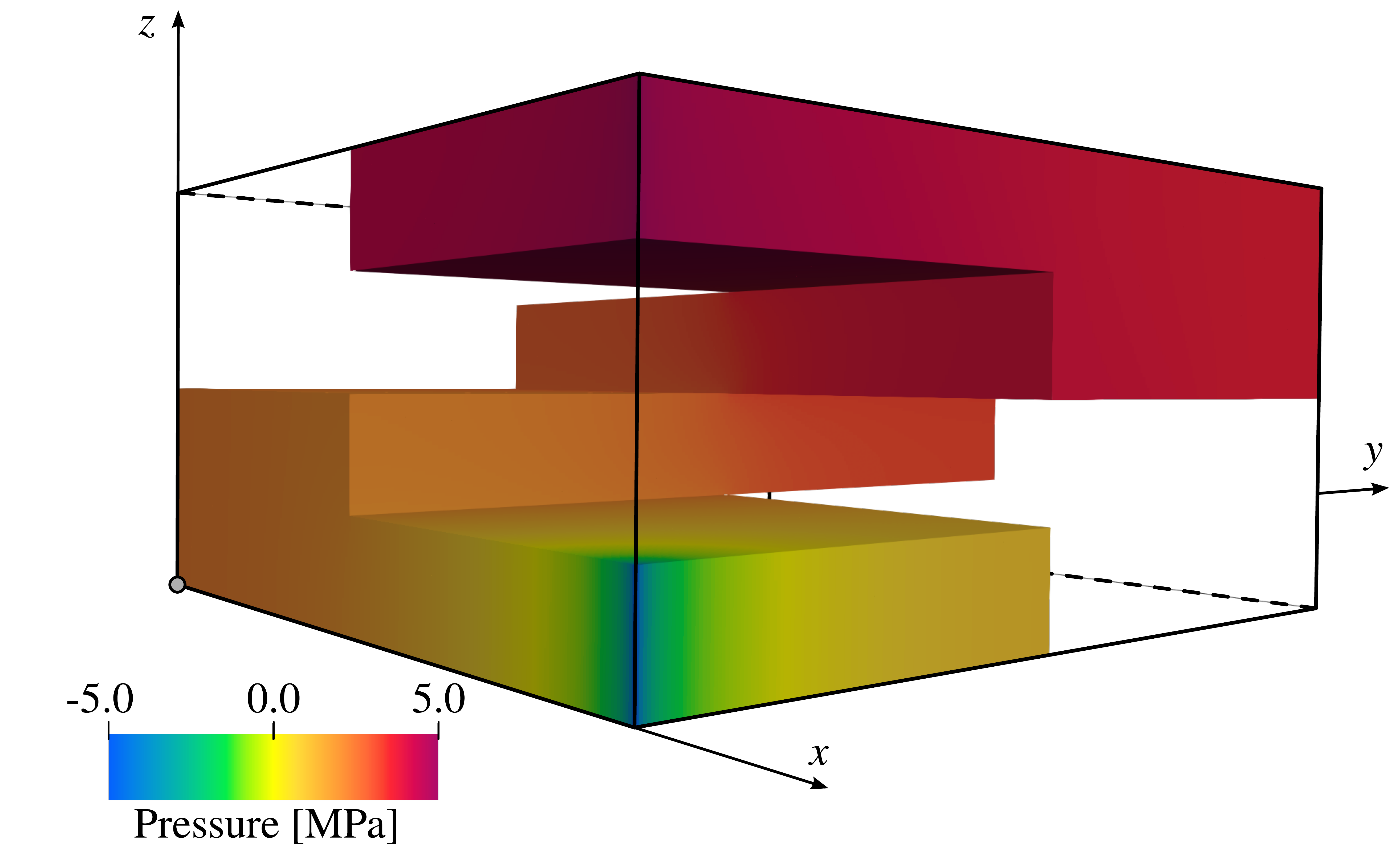}
        \caption{}
        \label{fig:staircase_pres}
    \end{subfigure}   
    \hfill\null
    \caption{Staircase benchmark configuration (a) and snapshots of the simulation results at $t=92$ days (b, c, d).   Note that in (a) injection and production wells are indicated in blue and red, respectively. }
    \label{fig:staircase_sketch}
\end{figure*}

\begin{table}[t]
        \caption{Weak scaling performance for the staircase example as a function of the number of cores and degrees of freedom (DoF) using different processor partitionings in (a) and (b).  Iteration counts are presented as average Newton iterations per timestep, and average GMRES iterations per Newton iteration.  Timings are broken into the preconditioner setup phase and the actual solve phase.  Note that setup time required for the preconditioner has two components: (i) a one-time expense needed to setup the elastic AMG preconditioner at the simulation start, and (ii) the recurring cost of setting up the flow preconditioner prior to each GMRES solve. Recurring costs are indicated as average times per Newton iteration.  The number of GMRES iterations per Newton iteration (in bold) exhibits only minor growth with mesh refinement, as desired.}

        \label{table:staircase}
        \centering
\small
(a)\\[5pt]
    \begin{tabular}{rrllcllrr}
    \toprule
    Cores & DoF & Iterations & && Setup Phase [s] && Solve Phase [s]\\
    \cline{3-4} \cline{6-7}
    && Newton per & GMRES per && Mechanics & Flow & \\
    && Timestep & Newton\\
   & & (avg.) & (avg.) && (once) & (avg.) & (avg).\\ 
    \midrule
    1 & 88,307 & 3.3 & \bf 13.5 && 0.39 & 0.05 & 0.38\\
    8 & 680,419 & 4.0 & \bf 14.4 && 1.12 &0.10 & 0.83\\
    64 & 5,342,147 &5.0 & \bf 16.9 && 1.51 & 0.18 & 1.85\\
    512 & 42,338,179 & 6.6 & \bf 21.1 && 2.66 & 0.23 & 2.60\\
    \bottomrule
    \end{tabular}\\[5pt]
    (b)\\[5pt]
        \begin{tabular}{rrllcllrr}
    \toprule
    Cores & DoF & Iterations & && Setup Phase [s] && Solve Phase [s]\\
    \cline{3-4} \cline{6-7}
    && Newton per & GMRES per && Mechanics & Flow & \\
    && Timestep & Newton\\
   & & (avg.) & (avg.) && (once) & (avg.) & (avg).\\ 
    \midrule
    2 & 88,307 & 3.3 & \bf 13.6 && 0.32 & 0.03 & 0.21\\
    16 & 680,419 & 4.0 & \bf 14.7 && 0.69 &0.06 & 0.42\\
    128 & 5,342,147 &5.0 & \bf 17.5 && 0.88 & 0.09 & 0.87\\
    1024 & 42,338,179 & 6.6 & \bf 21.7 && 1.75 & 0.16 & 1.30\\
    \bottomrule
    \end{tabular}
    
\end{table}

For the first benchmark, we consider the synthetic configuration illustrated in Figure~\ref{fig:staircase_sketch}.   Injection and production wells are placed at opposite ends of a high-permeability channel that winds its way in a ``staircase'' spiral through a low permeability host rock.  The reservoir channel has a 1000:1 permeability contrast with surrounding seal units, with isotropic permeability in both regions.  All external boundaries are assumed to be no flow.  For mechanics, the side and bottom boundaries have a zero normal-displacement constraint (roller boundary), while the top surface is allowed to deform freely.  The water injector and fluid production wells completely perforate the channel.  The injection (respectively, production) well reaches a target bottomhole pressure of 5 MPa above (below) initial formation pressure after 1 day of buildup (drawdown).  An illustrative snapshot of the simulation results is presented in Figure~\ref{fig:staircase_sketch}.  While the geometry is artificial, the problem is designed to exhibit strong coupling and provides a straightforward description for comparison with other codes.

The results of a weak-scaling study using this configuration are presented in Table~\ref{table:staircase}.  We begin with an  88k degree-of-freedom mesh run on one core.  The mesh is then uniformly refined and run on eight cores, preserving (approximately) the same number of unknowns per core.  This process is continued out to 512 cores and 42M degrees of freedom.  We then repeat the process with 88k degrees-of-freedom on 2 cores and scale out to 1024 cores.  One can also get an indication of strong scaling performance by comparing results in Tables~\ref{table:staircase}a and \ref{table:staircase}b for the same problem size but different core counts.

Because the problem is nonlinear, we observe an increase in the average number of Newton iterations per timestep as the mesh is refined.  As desired, however, the number of GMRES iterations per linear solve is relatively insensitive to the problem size, showing only minor growth with refinement.  Using more complex smoothers, we can further drive down the iteration count but at the cost of overall timing performance.  Table~\ref{table:staircase} also presents timing results.   Both the setup and solver iteration phases exhibit some growth, though even at large processor counts the solve time is quite satisfactory.  We note that the actual work per processor in this problem is quite small, and the timing growth is a symptom of communication costs in the AMG setup and solver iteration phases.  In general, however, the overall framework provides a robust basis for scalable performance.

\subsection{SPE10-based Benchmark}
\label{ex:spe10}

\begin{figure}[p]
    \centering
    \hfill
    \begin{subfigure}[t]{0.48\textwidth}
        \includegraphics[width=\textwidth]{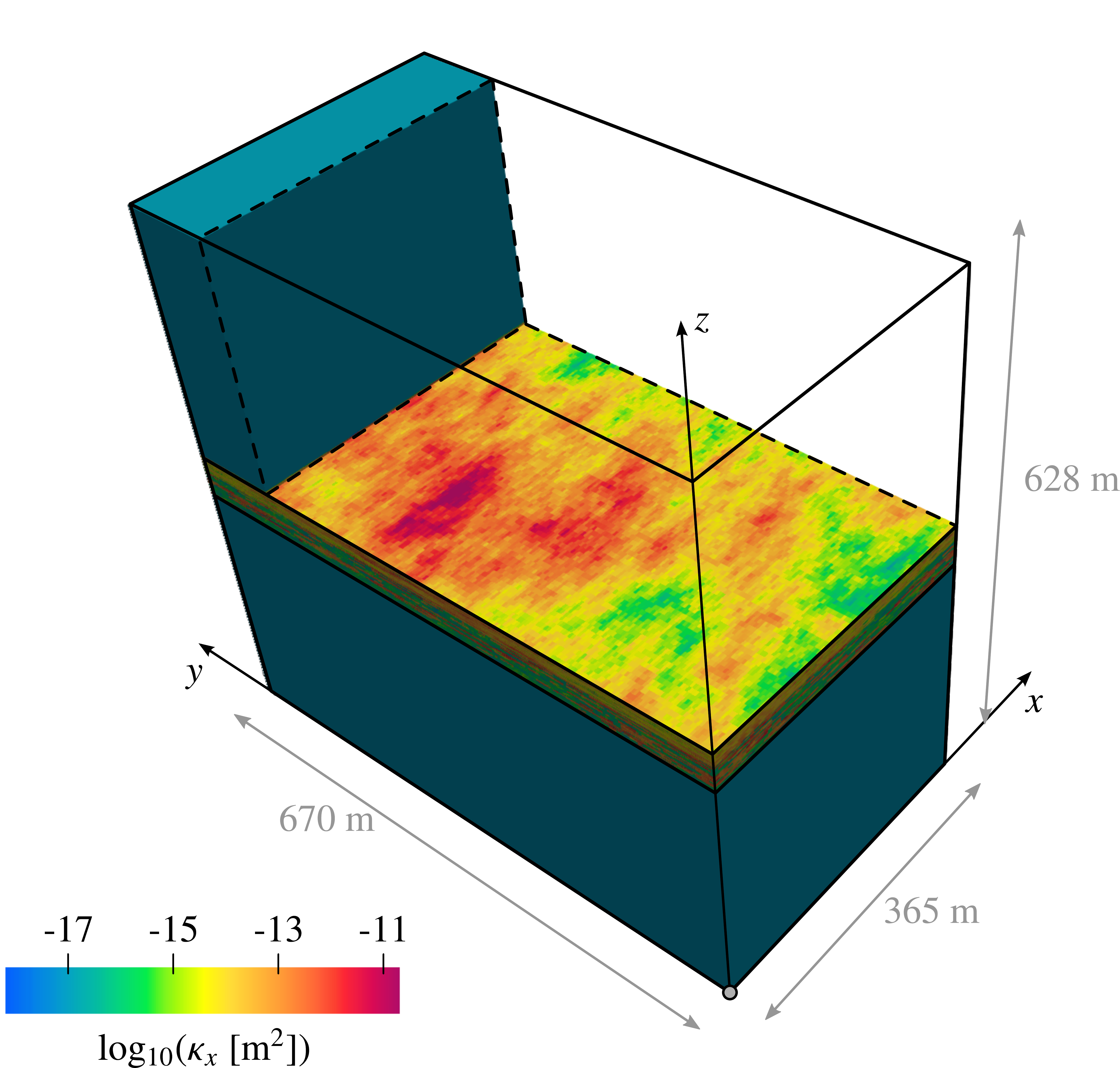}
        \caption{}
        \label{fig:spe10_permeability}
    \end{subfigure}
    \hfill
    \begin{subfigure}[t]{0.48\textwidth}
        \includegraphics[width=\textwidth]{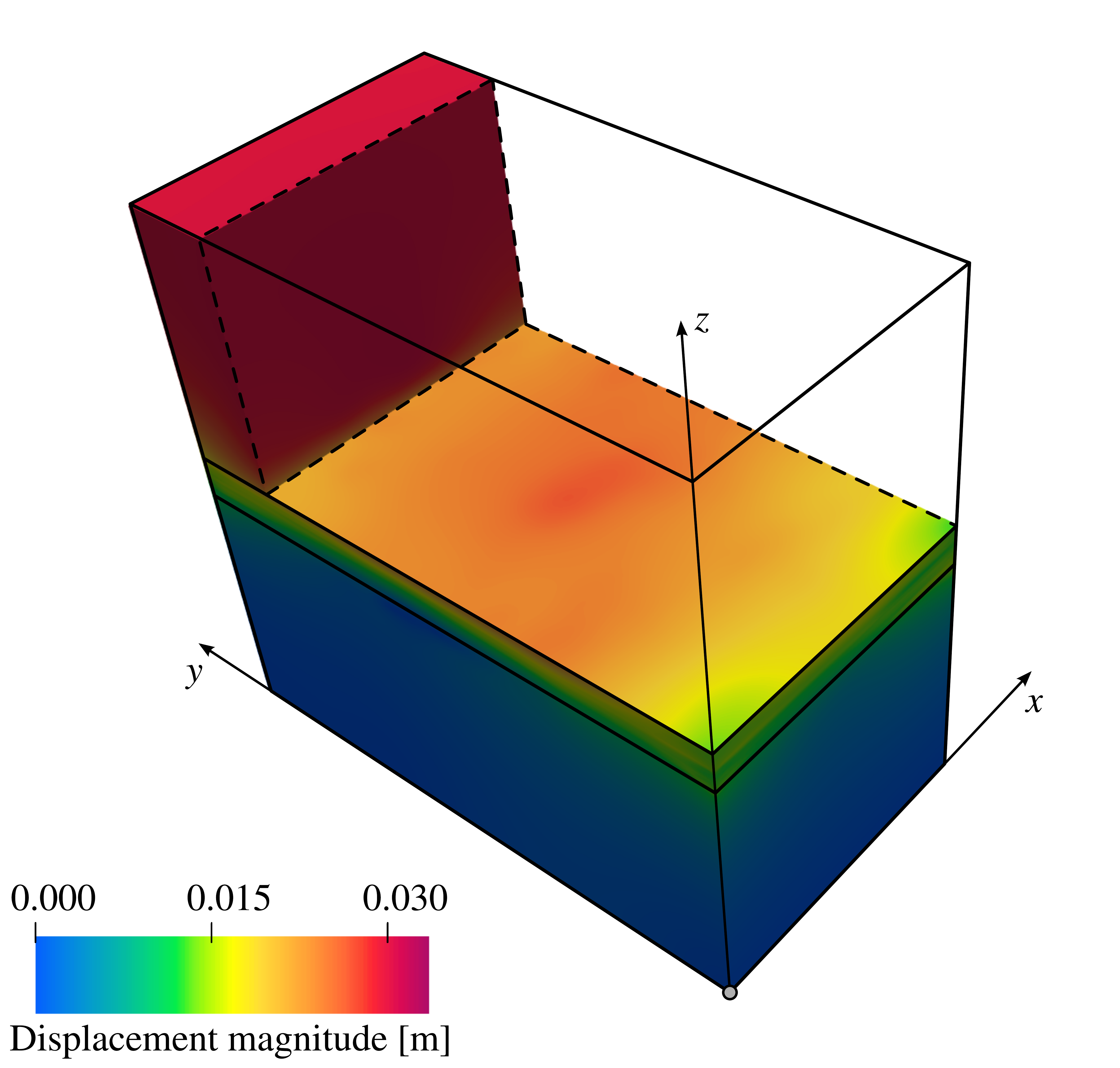}
        \caption{}
        \label{fig:spe10_disp}
    \end{subfigure}
    \hfill\null
    
    \hfill
    \begin{subfigure}[t]{0.48\textwidth}
        \includegraphics[width=\textwidth]{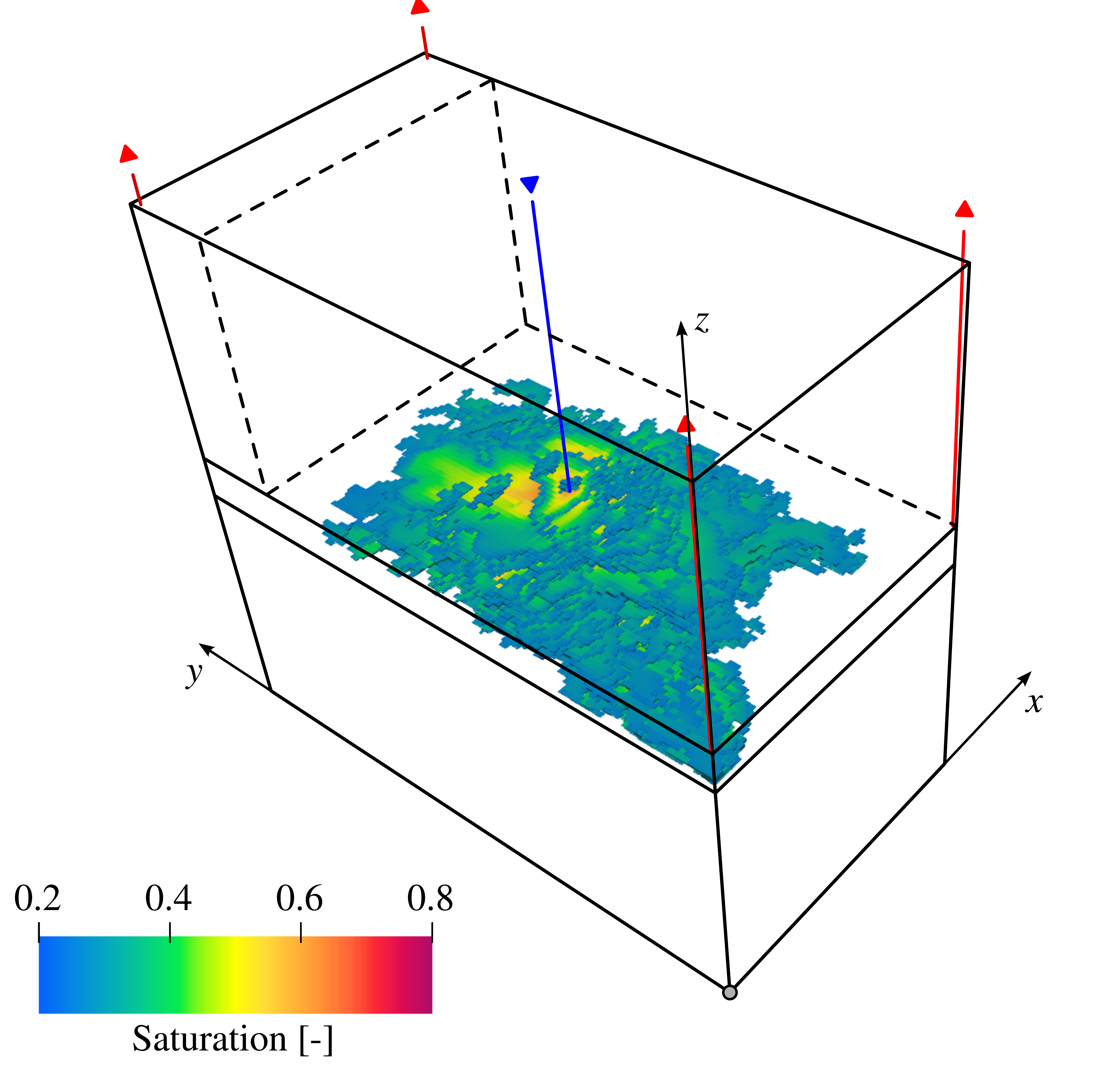}
        \caption{}
        \label{fig:spe10_sat}
    \end{subfigure}   
    \hfill
    \begin{subfigure}[t]{0.48\textwidth}
        \includegraphics[width=\textwidth]{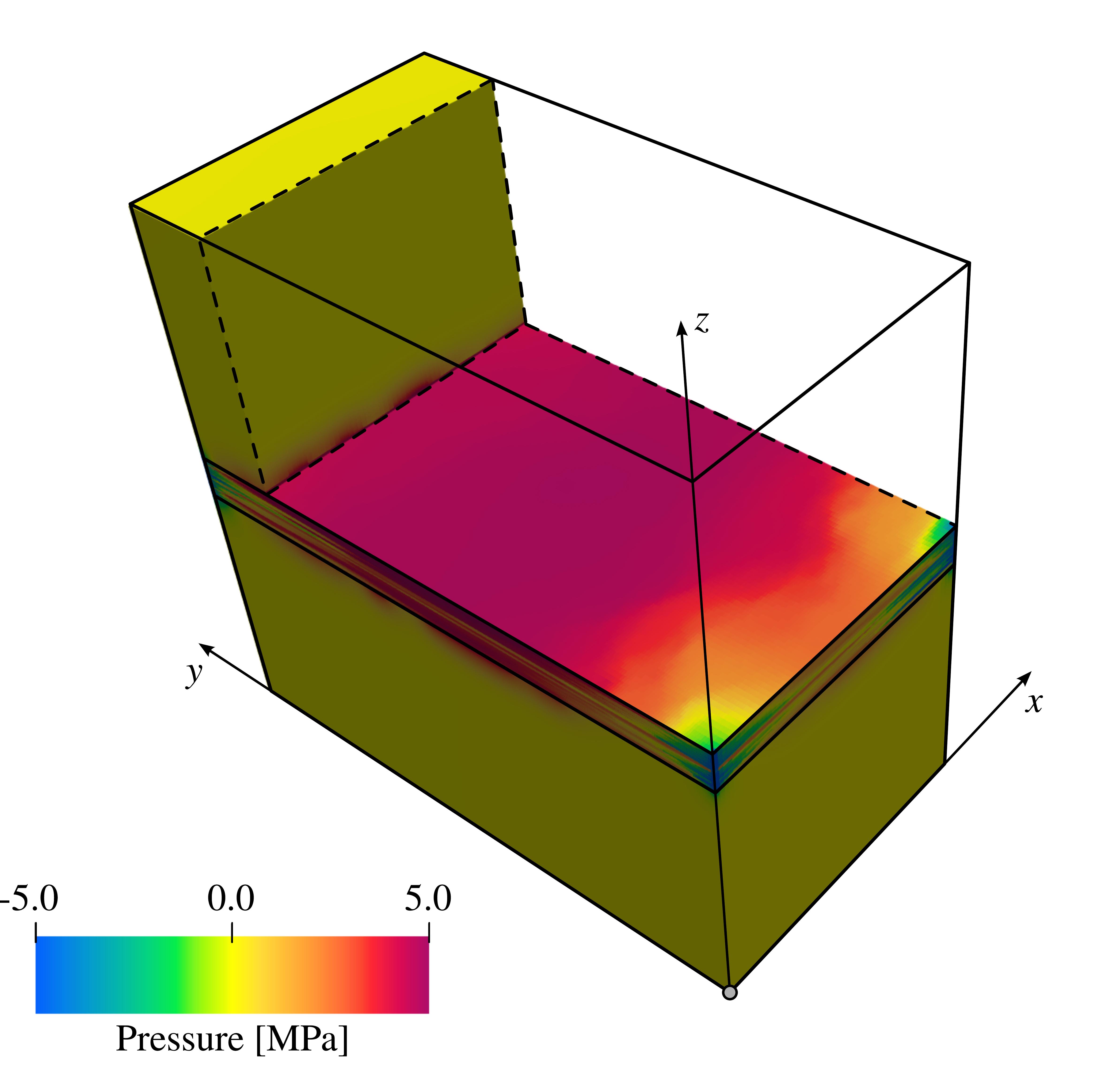}
        \caption{}
        \label{fig:spe10_pres}
    \end{subfigure}   
    \hfill\null
        \caption{Example 2: SPE10-based benchmark.  Sketch of the simulated domain and horizontal permeability field (a).   Snapshots of the simulation results at $t=98$ days (b, c, d).   Note that in (c) injection and production wells are indicated in blue and red, respectively.}
        \label{fig:spe10}
\end{figure}

As a more realistic test case, we consider a water flood in a heterogeneous reservoir.   The setup is based on Model~2 of the 10\textsuperscript{th} SPE Comparative Solution Project \cite{ChrBlu01} but equipped with poroelastic mechanical behavior.  In order to provide realistic mechanical boundary conditions, we add 300-m thick overburden and underburden layers to the original SPE10 reservoir.   A sketch of the simulated domain is shown in Figure \ref{fig:spe10}.   The computational grid consists of $60 \times 220 \times 252$ cells.  The total number of degrees of freedom is $n_{\text{total}}$ = 16,730,559 ($n_u$ = 10,077,759; $n_s$ = 3,326,400; and $n_p$ = 3,326,400).   The original SPE10 anisotropic permeability distribution and relative permeability curves are used---details provided in \cite{ChrBlu01}---while the porosity field has been slightly modified so that the lowest porosity for a cell is thresholded at 1\% to remove extremely low porosity values present in the original data files.   An isotropic permeability value equal to $0.01$ mD is assigned to the overburden and underburden layers, and porosity is set to 1\%.   Homogeneous Young's modulus $E$ = 5000 MPa, Poisson's ratio $\nu$ = 0.25, and Biot's coefficient $b$ = 1.0 are assumed everywhere.   All boundaries are impervious to fluid flow and constrained to have zero normal displacement, except for the top which is allowed to deform freely.   Five wells penetrate the entire reservoir as shown in Figure \ref{fig:spe10}.   The injector well, located in the middle of the domain, and the production wells, located at the corners, operate at a constant +5 and -5 MPa overpressure, respectively.  Note that these well conditions differ from the original SPE10 specification, which include rate-controlled wells.   Overall, 100 days of injection are simulated.   

\begin{table}[t]
        \caption{Strong scaling performance for the SPE10-based example.  The number of GMRES iterations is well-controlled and mostly insensitive to the parallel partitioning.}

        \label{tab:spe10}
\small
    \begin{tabular}{rrllcllrrr}
    \toprule
    Cores & DoF / Core & Iterations & && Setup Phase [s] && Solve Phase [s] & Overall [s] & Efficiency\\
    \cline{3-4} \cline{6-7}
    && Newton per & GMRES per && Mechanics & Flow & &\\
    && Timestep & Newton \\
   && (avg.)& (avg.)& & (once) & (avg.)  & (avg.)\\
    \midrule 
    
         36 & 464,738 & 6.5 & \bf 24.0 && 3.27 & 1.27	& 36.0 & 12,634 & 1.00 \\
     72 & 232,369 & 6.5 & \bf 24.7 &&  1.71 &  0.59 	& 18.0 & 6,303 & 1.00 \\
    144 & 116,184 & 6.5 & \bf 25.6 &&  0.99 &  0.31	& 9.0 &  3,149 & 1.00 \\
    288 &  58,092 & 6.5 & \bf 26.0 &&  1.36 &  0.18 	&  4.5 &  1,603 & 0.98 \\
    576 &  29,046 & 6.5 & \bf 28.1 &&  0.86 &  0.40 	&  2.2 &  881 & 0.90 \\
    
    \bottomrule
    \end{tabular}
\end{table}

Table \ref{tab:spe10} reports results of a strong-scaling study in which the total problem size is kept fixed but allocated across an increasing number of cores.  The GMRES iterations to convergence are well-controlled and mostly insensitive to the parallel partitioning.  Good overall timing efficiency is observed with only a mild decrease from ideal performance for the 576 core case.  The results suggest that one can use the proposed solver strategy and large processor counts to efficiently drive down long simulation times even for challenging heterogeneous problems.

\section{Conclusion}
\label{sec:conclusions}

In this work, we have presented a preconditioning scheme suitable for fully-implicit simulation of multiphase flow and transport with geomechanics.  We have used physics-based arguments to break the inherently coupled problem into modular sub-problems for which ``black-box'' algebraic methods can be effective.  This modularity also provides significant flexibility to test other approaches.  For example, while we have used a CPR-like strategy for the reduced multiphase flow problem, research into effective reservoir simulation preconditioners remains an active area of research.  If an alternative flow preconditioner has been identified, the partitioning framework here provides a simple strategy for adding poromechanical physics.

Regarding future work, a number of avenues would be interesting to explore.  First, it is well known that the coupling strength between different PDEs evolves with both time and timestep size.  Strong coupling is often observed at early time and small timesteps, while weaker coupling is observed at late time and large timesteps.  One might consider an adaptive solver strategy where the coupling strength between various quantities is detected based on the magnitude of coupling contributions in the preconditioner.  When sufficient ``uncoupling'' is detected between one or more fields via these coupling indicators, one might switch to a cheaper solver strategy.

Also, this work has focused almost exclusively on the solver strategy at the linearized equation level.  There are likely significant performance improvements to be gained by embedding similar ideas at the nonlinear level. Finally, while the physics-based splitting strategies used here are effective, they are PDE-specific and cannot be immediately extended to other coupled systems.  For example, the addition of thermal or fractured reservoir behavior would change the underlying PDEs and render the current solver strategy incomplete.  There is therefore significant room for both new physics-based solvers as well as algebraic-coupling methods that can be (to the extent possible) agnostic as to the underlying physics.

\section*{Acknowledgements}
Funding for JAW, NC, QMB and DO-K was provide by LLNL through Laboratory Directed Research and Development (LDRD) Project 18-ERD-027.  Funding for SK and HAT was provided by Total S.A. through the FC-MAELSTROM Project.   Portions of this work were performed under the auspices of the U.S. Department of Energy by Lawrence Livermore National Laboratory under Contract DE-AC52-07-NA27344.  


\appendix
\section{TPFA transmissibility computation}
\label{app:TPFA_computation}
\setcounter{figure}{0}  

\begin{figure}
    \centering
    \includegraphics[height=3.5cm]{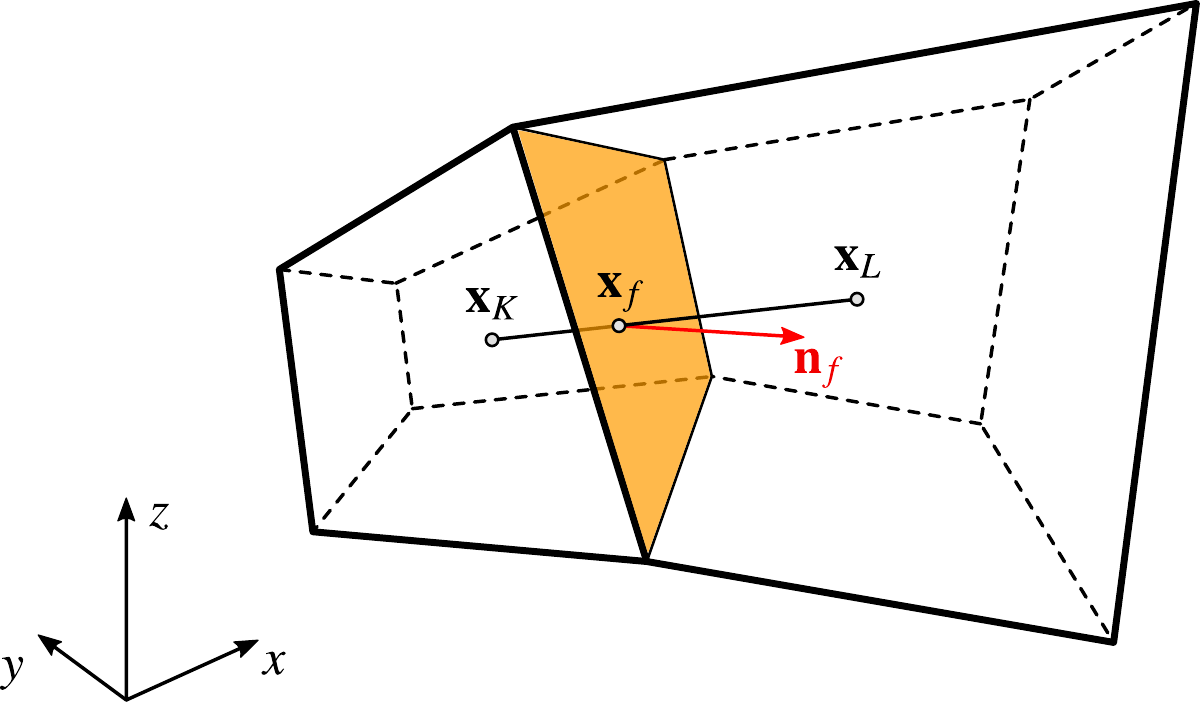}
    \caption{Geometrical entities involved in the TPFA transmissibility computation for an internal face shared by two cells.}
    \label{fig:TPFA}
\end{figure}

The transmissibility coefficient between two connected cells $K$ and $L$ in $\mathcal{T}^h$ sharing a face $f$ is defined by the following harmonic average

\begin{align}
   \Upsilon^f &= \frac{\Upsilon^{f,K}\Upsilon^{L,f}}{\Upsilon^{f,K} + \Upsilon^{L,f}},
   \label{eq:trans}
\end{align}

\noindent
where quantities $\Upsilon^{f,K}$ and $\Upsilon^{L,f}$, also known as half transmissibility coefficients, are given by \cite{TerVas13}

\begin{align}
   \Upsilon^{f,K} &= |f| \frac{(\vec{x}_f - \vec{x}_K) \cdot \tensorTwo{\kappa}_{\left| K \right.} \cdot \vec{n}_f}{(\vec{x}_f - \vec{x}_K) \cdot (\vec{x}_f - \vec{x}_K)}, &
   \Upsilon^{L,f} &= |f| \frac{(\vec{x}_L - \vec{x}_f) \cdot \tensorTwo{\kappa}_{\left| L \right.} \cdot \vec{n}_f}{(\vec{x}_L - \vec{x}_f) \cdot (\vec{x}_L - \vec{x}_2f)},
   \label{eq:half_trans}
\end{align}

\noindent
with $|f|$ the face area, and $\vec{x}_f$ a collocation point on $f$ that allows for enforcing point-wise pressure continuity across the face in case of incompressible single-phase flow.   In our implementation $\vec{x}_f$ is selected as the intersection of the face and the line connecting the cell centroids $\vec{x}_K$ and $\vec{x}_L$ \cite{KarDur16} as shown in Fig. \ref{fig:TPFA}.   For a boundary face $f \in \mathcal{F}_K \cap \Gamma_f^{D,h}$, transmissibility and half transmissibility coefficient coincide, i.e. $\Upsilon^f = \Upsilon^{f,K}$, and $\vec{x}_f$ is chosen as the orthogonal projection of $\vec{x}_K$ on $f$.   For a non-planar face, a plane characterized by a mean normal $\vec{n}_f$ and a collocation point $\vec{x}_f$ are associated to $f$.  

\section{Finite element and finite volume vectors and matrices}
\label{app:FEFV_vec_mat_blocks}
Repeated indices in the expressions given in this appendix do not imply summation.

\subsection{Discrete residual vectors expressions}
\label{app:discr_res}
The residual vector $\Vec{r}_u$~\eqref{eq:NL_res} is assembled as sum of element contributions.   The residual vector $\Vec{r}_\ell$~\eqref{eq:NL_res}, $\ell = \{ w, nw \}$, is assembled as sum of element ($\Vec{r}^K_\ell$) and interface ($\Vec{r}^f_\ell$) contributions, i.e. $\Vec{r}_\ell = \Vec{r}^K_\ell + \Vec{r}^f_\ell$.   Their expressions read:
\begin{align}
  \left[ \Vec{r}_n^{u,K} \right]_i =&\; 
	\int_K \nabla^s \tensorTwo{\eta}_i : \tensorTwo{\sigma}^\prime_n \; \mathrm{d}V -
	\int_K \nabla \cdot \tensorTwo{\eta}_i b p_n \; \mathrm{d}V -
	\int_K \tensorTwo{\eta}_i \cdot \rho_n \vec{g} \; \mathrm{d}V -
	\int_{\partial K \cap \Gamma_u^{N,h}} \tensorTwo{\eta}_i \cdot \bar{\vec{t}}_n \; \mathrm{d}A
	&&\forall i \in \{1, 2, \ldots, n_{u} \} \\
  \left[ \Vec{r}_n^{w,K}  \right]_i =&\; 
  \int_K \psi_i \left[ m_{w, n} - m_{w, (n-1)}\right] \; \mathrm{d}V -
  \Delta t_n \int_K \psi_i \left( q_{w,n}^I - q_{w,n}^P \right) \, \mathrm{d}V
  &&\forall i \in \{1, 2, \ldots, n_s \} \\
  \left[ \Vec{r}_n^{w,f}\right]_i =&\;
  - \Delta t_n \sum_{f \notin \Gamma_f^{N,h}} \llbracket \psi_i \rrbracket_f  F_{w,n}^f
  &&\forall i \in \{1, 2, \ldots, n_s \} \\
  \left[ \Vec{r}_n^{nw,K}  \right]_i =&\; 
  \int_K \chi_i \left[ m_{nw, n} - m_{nw, (n-1)}\right] \; \mathrm{d}V -
  \Delta t_n \int_K \chi_i \left( q_{nw,n}^I - q_{nw,n}^P \right) \, \mathrm{d}V
  &&\forall i \in \{1, 2, \ldots, n_p \} \\
  \left[ \Vec{r}_n^{nw,f}\right]_i =&\;
  - \Delta t_n \sum_{f \notin \Gamma_f^{N,h}} \llbracket \chi_i \rrbracket_f  F_{nw,n}^f
  &&\forall i \in \{1, 2, \ldots, n_p\}
\end{align}

\subsection{Jacobian sub-matrices expressions}
The sub-matrices appearing in the global Jacobian matrix~\eqref{eq:jac_system} are constructed as sum of element ($\Mat{A}^K_{lm}$) and face ($\Mat{A}^f_{lm}$) contributions, namely
\begin{align}
	\Mat{A}_{lm} &= \Mat{A}_{lm}^K + \Mat{A}_{lm}^f, &&\forall (l,m) \in \{u,s,p\} \times \{u,s,p\}.
	\label{eq:jac_assembly}
\end{align}

\noindent
The expressions for sub-matrices $\Mat{A}_{lm}^K$ read:
\begin{subequations}
\begin{align}
  \left[\Mat{A}_{uu}^K\right]_{i,j} =&\; 
  \int_K \nabla^s \tensorTwo{\eta}_i : \tensorFour{C}_{dr} : \nabla^s \tensorTwo{\eta}_j \, \mathrm{d}V
  -\int_K \tensorTwo{\eta}_i \cdot \left(\left.\frac{\partial \rho}{\partial \epsilon_v}\right|_n \vec{g} \right) \nabla \cdot \tensorTwo{\eta}_j \, \mathrm{d}V 
  &&\forall(i,j) \in \{1, 2, \ldots, n_{u} \} \times \{1, 2, \ldots, n_{u} \}\\
  \left[\Mat{A}_{us}^K\right]_{i,j} =&\; 
  -\int_K \tensorTwo{\eta}_i \cdot \left(\left.\frac{\partial \rho}{\partial s}\right|_n \vec{g} \right) \psi_j \, \mathrm{d}V
  &&\forall(i,j) \in \{1, 2, \ldots, n_{u} \} \times \{1, 2, \ldots, n_s \}\\
  \left[\Mat{A}_{up}^K\right]_{i,j} =&\;
  -\int_K \nabla \cdot \tensorTwo{\eta}_i  b \chi_j \, \mathrm{d}V
  -\int_K \tensorTwo{\eta}_i \cdot \left(\left.\frac{\partial \rho}{\partial p}\right|_n \vec{g} \right) \chi_j \, \mathrm{d}V 
  &&\forall(i,j) \in \{1, 2, \ldots, n_{u} \} \times \{1, 2, \ldots, n_p \}\\
  \left[\Mat{A}_{su}^K\right]_{i,j} =&\; 
  \int_K \psi_i \left( \left.\frac{\partial m_{w}}{\partial \epsilon_v}\right|_n \right) \nabla \cdot \tensorTwo{\eta}_j  \, \mathrm{d}V
  &&\forall(i,j) \in \{1, 2, \ldots, n_s \} \times \{1, 2, \ldots, n_{u} \}\\
  \left[\Mat{A}_{ss}^K\right]_{i,j} =&\;
  \int_K \psi_i \left( \left.\frac{\partial m_{w}}{\partial s}\right|_n \right) \psi_j \, \mathrm{d}V
  -\Delta t \int_K \psi_i \left( \left.\frac{\partial q^I_{w}}{\partial s}\right|_n - \left.\frac{\partial q^P_{w}}{\partial s}\right|_n\right) \psi_j \, \mathrm{d}V
  &&\forall(i,j) \in \{1, 2, \ldots, n_s \} \times \{1, 2, \ldots, n_s \} \\
  \left[\Mat{A}_{sp}^K\right]_{i,j} =&\;
  \int_K \psi_i \left( \left.\frac{\partial m_{w}}{\partial p}\right|_n \right) \chi_j \, \mathrm{d}V
  -\Delta t \int_K \psi_i \left( \left.\frac{\partial q^I_{w}}{\partial p}\right|_n - \left.\frac{\partial q^P_{w}}{\partial p}\right|_n \right) \chi_j \, \mathrm{d}V
  &&\forall(i,j) \in \{1, 2, \ldots, n_s \} \times \{1, 2, \ldots, n_p \}\\ 
  \left[\Mat{A}_{pu}^K\right]_{i,j} =&\; 
  \int_K \chi_i \left( \left.\frac{\partial m_{nw}}{\partial \epsilon_v}\right|_n \right) \nabla \cdot \tensorTwo{\eta}_j  \, \mathrm{d}V
  &&\forall(i,j) \in \{1, 2, \ldots, n_p \} \times \{1, 2, \ldots, n_u \}\\
  \left[\Mat{A}_{ps}^K\right]_{i,j} =&\;
  \int_K \chi_i \left( \left.\frac{\partial m_{nw}}{\partial s}\right|_n \right) \psi_j \, \mathrm{d}V
  -\Delta t \int_K \chi_i \left( \left.\frac{\partial q^I_{nw}}{\partial s}\right|_n - \left.\frac{\partial q^P_{nw}}{\partial s}\right|_n \right) \psi_j \, \mathrm{d}V
  &&\forall(i,j) \in \{1, 2, \ldots, n_p \} \times \{1, 2, \ldots, n_s \}\\  
  \left[\Mat{A}_{pp}^K\right]_{i,j} =&\;
  \int_K \chi_i \left( \left.\frac{\partial m_{nw}}{\partial p}\right|_n \right) \chi_j \, \mathrm{d}V
  -\Delta t \int_K \chi_i \left( \left.\frac{\partial q^I_{nw}}{\partial p}\right|_n - \left.\frac{\partial q^P_{nw}}{\partial p}\right|_n \right) \chi_j \, \mathrm{d}V
  &&\forall(i,j) \in \{1, 2, \ldots, n_p \} \times \{1, 2, \ldots, n_p \}                          
\end{align}
\end{subequations}

\noindent
The partial derivatives in the above equations are expanded as:

\begin{subequations}
\begin{align}
  \left.\frac{\partial \rho}{\partial \epsilon_v}\right|_n =&\; \left[-\rho_{s,n} + s_n \rho_{w,n} + \left( 1 - s_n \right) \rho_{nw,n} \right] \left.\frac{\partial \phi}{\partial \epsilon_v}\right|_n  \\
  \left.\frac{\partial \rho}{\partial s}\right|_n =&\; \phi_n \left(\rho_{w,n} - \rho_{nw,n} \right) \\
  \left.\frac{\partial \rho}{\partial p}\right|_n =&\; \left[-\rho_{s,n} + s_n \rho_{w,n} + \left( 1 - s_n \right) \rho_{nw,n} \right] \left.\frac{\partial \phi}{\partial p}\right|_n +
  (1-\phi_n) \left.\frac{\partial \rho_s}{\partial p}\right|_n + s_n \phi_n  \left.\frac{\partial \rho_w}{\partial p}\right|_n + (1-s_n) \phi_n \left.\frac{\partial \rho_{nw}}{\partial p}\right|_n \\
  \left.\frac{\partial m_{w}}{\partial \epsilon_v}\right|_n =&\; s_n \rho_{w,n} \left.\frac{\partial \phi}{\partial \epsilon_v}\right|_n \\ 
  \left.\frac{\partial m_{w}}{\partial s}\right|_n =&\; \phi_n \rho_{w,n} \\
  \left.\frac{\partial q^I_{w}}{\partial s}\right|_n =&\; \rho_{w,n} \sum_{\ell=\{w,nw\}} \left( \frac{1}{\mu_{\ell,n}} \left.\frac{\partial k_{r\ell}}{\partial s}\right|_n \right) \Phi^{K,W}_{w,n} \; \delta \left(\vec{x} - \vec{x}_K \right)\\
  \left.\frac{\partial q^P_{w}}{\partial s}\right|_n =&\; - \rho_{w,n} \frac{1}{\mu_{w,n}} \left.\frac{\partial k_{rw}}{\partial s}\right|_n  \Phi^{K,W}_{w,n} \; \delta \left(\vec{x} - \vec{x}_K \right)\\
  \left.\frac{\partial m_{w}}{\partial p}\right|_n =&\; s_n \rho_{w,n} \left.\frac{\partial \phi}{\partial p}\right|_n + \phi_n s_n \left.\frac{\partial \rho_w}{\partial p}\right|_n \\
  \left.\frac{\partial q^I_{w}}{\partial p}\right|_n =&\;
    \left.\frac{\partial \rho_{w}}{\partial p}\right|_n \sum_{\ell=\{w,nw\}} \left( \lambda_{\ell,n} \right) \Phi^{K,W}_{w,n} \; \delta \left(\vec{x} - \vec{x}_K \right) -
    \rho_{w,n} \sum_{\ell=\{w,nw\}} \left( \frac{\lambda_{\ell,n}}{\mu_{\ell,n}} \left.\frac{\partial \mu_{\ell}}{\partial p}\right|_n \right) \Phi^{K,W}_{w,n} \; \delta \left(\vec{x} - \vec{x}_K \right)
    \nonumber \\
    &\; + \rho_{w,n} WI \sum_{\ell=\{w,nw\}} \left( \lambda_{\ell,n} \right) \left[-1 + \left.\frac{\partial \rho_{w}}{\partial p}\right|_n g  \left( z_{bh} - z \right)  \right] \; \delta \left(\vec{x} - \vec{x}_K \right) \\
  \left.\frac{\partial q^P_{w}}{\partial p}\right|_n =&\; 
  -\left.\frac{\partial \rho_{w}}{\partial p}\right|_n \lambda_{w,n} \Phi^{K,W}_{w,n} \; \delta \left(\vec{x} - \vec{x}_K \right)
  +\rho_{w,n} \frac{\lambda_{w,n}}{\mu_{w,n}} \left.\frac{\partial \mu_{w}}{\partial p}\right|_n \Phi^{K,W}_{w,n} \; \delta \left(\vec{x} - \vec{x}_K \right) 
  \nonumber \\
  &\; - \rho_{w,n} \lambda_{w,n} WI \left[-1 + \left.\frac{\partial \rho_{w}}{\partial p}\right|_n g  \left( z_{bh} - z \right)  \right] \; \delta \left(\vec{x} - \vec{x}_K \right) \\
  \left.\frac{\partial m_{nw}}{\partial \epsilon_v}\right|_n =&\; \left(1 - s_n\right) \rho_{nw,n} \left.\frac{\partial \phi}{\partial \epsilon_v}\right|_n \\  
  \left.\frac{\partial m_{nw}}{\partial s}\right|_n =&\; -\phi_n \rho_{nw,n} \\
  \left.\frac{\partial q^I_{nw}}{\partial s}\right|_n =&\; \rho_{nw,n} \sum_{\ell=\{w,nw\}} \left( \frac{1}{\mu_{\ell,n}} \left.\frac{\partial k_{r\ell}}{\partial s}\right|_n \right) \Phi^{K,W}_{nw,n} \; \delta \left(\vec{x} - \vec{x}_K \right)\\
  \left.\frac{\partial q^P_{nw}}{\partial s}\right|_n =&\; - \rho_{nw,n} \frac{1}{\mu_{nw,n}} \left.\frac{\partial k_{rnw}}{\partial s}\right|_n  \Phi^{K,W}_{nw,n} \; \delta \left(\vec{x} - \vec{x}_K \right)\\
  \left.\frac{\partial m_{nw}}{\partial          p}\right|_n =&\; \left(1 - s_n\right) \rho_{nw,n} \left.\frac{\partial \phi}{\partial p}\right|_n + \phi_n \left(1 - s_n\right) \left.\frac{\partial \rho_{nw}}{\partial p}\right|_n \\
  \left.\frac{\partial q^I_{nw}}{\partial p}\right|_n =&\;
  \left.\frac{\partial \rho_{nw}}{\partial p}\right|_n \sum_{\ell=\{w,nw\}} \left( \lambda_{\ell,n} \right) \Phi^{K,W}_{nw,n} \; \delta \left(\vec{x} - \vec{x}_K \right)-
  \rho_{nw,n} \sum_{\ell=\{w,nw\}} \left( \frac{\lambda_{\ell,n}}{\mu_{\ell,n}} \left.\frac{\partial \mu_{\ell}}{\partial p}\right|_n \right) \Phi^{K,W}_{nw,n}  \; \delta \left(\vec{x} - \vec{x}_K \right)
  \nonumber \\
  &\;  
  + \rho_{nw,n} WI \sum_{\ell=\{w,nw\}} \left( \lambda_{\ell,n} \right) \left[-1 + \left.\frac{\partial \rho_{nw}}{\partial p}\right|_n g  \left( z_{bh} - z \right)  \right]  \; \delta \left(\vec{x} - \vec{x}_K \right) \\
  \left.\frac{\partial q^P_{nw}}{\partial p}\right|_n =&\; 
  -\left.\frac{\partial \rho_{nw}}{\partial p}\right|_n \lambda_{nw,n} \Phi^{K,W}_{nw,n} \; \delta \left(\vec{x} - \vec{x}_K \right) +
  \rho_{nw,n} \frac{\lambda_{nw,n}}{\mu_{nw,n}} \left.\frac{\partial \mu_{nw}}{\partial p}\right|_n \Phi^{K,W}_{nw,n} \; \delta \left(\vec{x} - \vec{x}_K \right) 
  \nonumber \\
  &\;  
  - \rho_{nw,n} \lambda_{nw,n} WI \left[-1 + \left.\frac{\partial \rho_{nw}}{\partial p}\right|_n g  \left( z_{bh} - z \right)  \right]  \; \delta \left(\vec{x} - \vec{x}_K \right)
\end{align}
\end{subequations}

\noindent
The finite volume scheme provides a direct functional dependence for intercell fluxes $F^f_{\ell,n} $ on the algebraic vectors $ \Vec{s_n}$ and  $\Vec{p_n}$. The global expressions for sub-matrices $\Mat{A}_{lk}^f$ are therefore\\

\begin{subequations}
\begin{align}
  \left[\Mat{A}_{uu}^f\right]_{i,j} =&\; 0 
  &&\forall(i,j) \in \{1, 2, \ldots, n_{u} \} \times \{1, 2, \ldots, n_{u} \}\\
  \left[\Mat{A}_{us}^f\right]_{i,j} =&\; 0
  &&\forall(i,j) \in \{1, 2, \ldots, n_{u} \} \times \{1, 2, \ldots, n_s \}\\
  \left[\Mat{A}_{up}^f\right]_{i,j} =&\; 0
  &&\forall(i,j) \in \{1, 2, \ldots, n_{u} \} \times \{1, 2, \ldots, n_p \}\\
  \left[\Mat{A}_{su}^f\right]_{i,j} =&\; 0
  &&\forall(i,j) \in \{1, 2, \ldots, n_s \} \times \{1, 2, \ldots, n_{u} \} \\
  \left[\Mat{A}_{ss}^f\right]_{i,j} =&\; \llbracket \psi_i \rrbracket_f \Delta t_n \left.\frac{\partial F_{w}^f}{\partial {s}_j}\right|_n
  &&\forall(i,j) \in \{1, 2, \ldots, n_s \} \times \{1, 2, \ldots, n_s \} \\
  \left[\Mat{A}_{sp}^f\right]_{i,j} =&\; \llbracket \psi_i \rrbracket_f \Delta t_n \left.\frac{\partial F_{w}^f}{\partial {p}_j}\right|_n
  &&\forall(i,j) \in \{1, 2, \ldots, n_s \} \times \{1, 2, \ldots, n_p \} \\
  \left[\Mat{A}_{pu}^f\right]_{i,j} =&\; 0
  &&\forall(i,j) \in \{1, 2, \ldots, n_p \} \times \{1, 2, \ldots, n_{u} \} \\
  \left[\Mat{A}_{ps}^f\right]_{i,j} =&\; \llbracket \chi_i \rrbracket_f \Delta t_n \left.\frac{\partial F_{nw}^f}{\partial {s}_j}\right|_n
  &&\forall(i,j) \in \{1, 2, \ldots, n_p \} \times \{1, 2, \ldots, n_s \} \\
  \left[\Mat{A}_{pp}^f\right]_{i,j} =&\; \llbracket \chi_i \rrbracket_f \Delta t_n \left.\frac{\partial F_{nw}^f}{\partial {p}_j}\right|_n
  &&\forall(i,j) \in \{1, 2, \ldots, n_p \} \times \{1, 2, \ldots, n_p \}
\end{align}
\label{eq:jac_f_int}
\end{subequations}

\noindent
The partial derivatives in the above equations are expanded as:

\begin{subequations}
\begin{align}
  \left.\frac{\partial F_{w}^f}{\partial {s}_j}\right|_n &= \left( \frac{\rho^{\texttt{upw}}_{w,n}}{\mu^{\texttt{upw}}_{w,n}} \Phi^{K,f}_{w,n} \right) \left.\frac{\partial k^{\texttt{upw}}_{rw}}{\partial {s}_j}\right|_n \\
  \left.\frac{\partial F_{w}^f}{\partial {p}_j}\right|_n &= \left( \left. \frac{\partial \rho^{\texttt{upw}}_{w}}{\partial {p}_j} \right|_n \lambda^{\texttt{upw}}_{w,n} - \rho^{\texttt{upw}}_{w,n} \frac{\lambda^{\texttt{upw}}_{w,n}}{\mu^{\texttt{upw}}_{w,n}} \left. \frac{\partial \mu^{\texttt{upw}}_{w}}{\partial {p}_j} \right|_n \right) \Phi^f_{w,n} + \rho^{\texttt{upw}}_{w,n} \lambda^{\texttt{upw}}_{w,n} \Upsilon^f \left( \left. \frac{\partial p_L}{\partial {p}_j} \right|_n - \left. \frac{\partial p_K}{\partial {p}_j} \right|_n + \left. \frac{\partial \varrho^f_{w}}{\partial {p}_j} \right|_n g (z_L - z_K) \right) \\
  \left.\frac{\partial F_{nw}^f}{\partial {s}_j}\right|_n &= \left( \frac{\rho^{\texttt{upw}}_{nw,n}}{\mu^{\texttt{upw}}_{nw,n}} \Phi^{K,f}_{nw,n} \right) \left.\frac{\partial k^{\texttt{upw}}_{rnw}}{\partial {s}_j}\right|_n \\
  \left.\frac{\partial F_{nw}^f}{\partial {p}_j}\right|_n &= \left( \left. \frac{\partial \rho^{\texttt{upw}}_{nw}}{\partial {p}_j} \right|_n \lambda^{\texttt{upw}}_{nw,n} - \rho^{\texttt{upw}}_{nw,n} \frac{\lambda^{\texttt{upw}}_{nw,n}}{\mu^{\texttt{upw}}_{nw,n}} \left. \frac{\partial \mu^{\texttt{upw}}_{nw}}{\partial {p}_j} \right|_n \right) \Phi^f_{nw,n} + \rho^{\texttt{upw}}_{nw,n} \lambda^{\texttt{upw}}_{nw,n} \Upsilon^f \left( \left. \frac{\partial p_L}{\partial {p}_j} \right|_n - \left. \frac{\partial p_K}{\partial {p}_j} \right|_n + \left. \frac{\partial \varrho^f_{nw}}{\partial {p}_j} \right|_n g (z_L - z_K) \right)
\end{align}
\end{subequations}




\bibliography{biblio_CPR_fixed_stress}



\end{document}